\documentclass[leqno,a4paper,12pt]{article}

\usepackage{amsmath,amsthm}

\newtheorem{Thm}{\indent Theorem}[section]
\newtheorem{Prop}[Thm]{\indent Proposition}
\newtheorem{Lemma}[Thm]{\indent Lemma}
\newtheorem{Cor}[Thm]{\indent Corollary}
\newtheorem{Comp}[Thm]{\indent Complement}
\newtheorem{Var}[Thm]{\indent Variant}

\theoremstyle{definition}
\newtheorem{Rem}[Thm]{\indent Remark}

\newtheorem{Def}[Thm]{\indent Definition}
\newtheorem{Ex}[Thm]{\indent Example}

\def\qed{{\hskip0pt\unskip\unskip\nobreak\hfil\penalty50
          \hskip1em\hbox{}\nobreak\hfil
          {\bf q.e.d.}%
          \parfillskip=0pt\finalhyphendemerits=0
          \par}\medskip}

\newenvironment{Proof}
               {{\it Proof.}\quad}
               {\qed}

\newenvironment{Proofof}[1]
               {{\it Proof of #1.}\quad}
               {\qed}

\newcommand{\Prime}{\kern3\fontdimen1\font$'$\kern-7\fontdimen1\font}

\long\def\forget#1{}

\long\def\beginSIDEREMARK#1\endSIDEREMARK
    {{\par\bigskip\advance\leftskip by 2cm
                  \advance\rightskip by -2cm\noindent
      {\bf Our own side remark:} #1
      \par\bigskip\noindent}}

\long\def\beginFORGET#1\endFORGET{#1}
\long\def\beginFORGET#1\endFORGET{}

\def\?{\ ???\ \immediate\write16{}%
\immediate\write16{Warning: There was still a question mark . . . }%
\immediate\write16{}}

\usepackage{amsmath}

\usepackage{exscale}
\usepackage{amssymb}

\newcommand{\BA}{{\mathbb{A}}}

\newcommand{\BC}{{\mathbb{C}}}

\newcommand{\BG}{{\mathbb{G}}}

\newcommand{\BQ}{{\mathbb{Q}}}
\newcommand{\BR}{{\mathbb{R}}}

\newcommand{\BZ}{{\mathbb{Z}}}

\newcommand{\FF}{{\mathfrak{F}}}

\newcommand{\FH}{{\mathfrak{H}}}

\newcommand{\FS}{{\mathfrak{S}}}
\newcommand{\FT}{{\mathfrak{T}}}
\newcommand{\FU}{{\mathfrak{U}}}
\newcommand{\FV}{{\mathfrak{V}}}
\newcommand{\FW}{{\mathfrak{W}}}
\newcommand{\FX}{{\mathfrak{X}}}
\newcommand{\FY}{{\mathfrak{Y}}}
\newcommand{\FZ}{{\mathfrak{Z}}}

\newcommand{\CA}{{\cal A}}
\newcommand{\CB}{{\cal B}}
\newcommand{\CC}{{\cal C}}

\newcommand{\Stab}{\mathop{\rm Stab}\nolimits}
\newcommand{\Hom}{\mathop{\rm Hom}\nolimits}
\newcommand{\kernel}{\mathop{\rm Ker}\nolimits}
\newcommand{\loccit}{[loc.$\;$cit.]}
\newcommand{\Res}{\mathop{\rm Res}\nolimits}

\def\tei{\, | \,}
\def\halb{\frac{1}{2}}

\def\id{{\rm id}}

\def\invlim{\mathop{\vtop{\hbox{\rm lim}\vskip-8pt
        \hbox{\hskip1pt$\scriptstyle\longleftarrow$}\vskip-1pt}}}

\newbox\mybox
\def\arrover#1{\mathrel{
       \setbox\mybox=\hbox spread 1.4em{\hfil$\scriptstyle#1$\hfil}
       \vbox{\offinterlineskip\copy\mybox
             \hbox to\wd\mybox{\rightarrowfill}}}}
\def\larrover#1{\mathrel{
       \setbox\mybox=\hbox spread 1.4em{\hfil$\scriptstyle#1$\hfil}
       \vbox{\offinterlineskip\copy\mybox
             \hbox to\wd\mybox{\leftarrowfill}}}}

\def\ontoover#1{\mathrel{
       \setbox\mybox=\hbox spread 1.4em{\hfil$\scriptstyle#1$\hfil}
       \vbox{\offinterlineskip\copy\mybox
             \hbox to\wd\mybox{\rightarrowfill\hskip-2.8mm
                               $\rightarrow$}}}}
\def\leftontoover#1{\mathrel{
       \setbox\mybox=\hbox spread 1.4em{\hfil$\scriptstyle#1$\hfil}
       \vbox{\offinterlineskip\copy\mybox
             \hbox to\wd\mybox{$\leftarrow$\hskip-2.8mm
                               \leftarrowfill}}}}
\def\longto{\longrightarrow}
\def\into{\hookrightarrow}

\def\longonto{\ontoover{\ }}
\def\isoto{\arrover{\sim}}

\def\longinto{\lhook\joinrel\longrightarrow}

\usepackage[curve,matrix,arrow,cmtip]{xy}
\NoComputerModernTips

\def\myxymessage{\def\messagetext
   {Here an xy-pic diagram was omitted to speed up compilation . . . }
   \immediate\write16{\messagetext}
   \hbox{\bf \messagetext}}
\def\filxymatrix#1{\myxymessage}
\def\filxyarray#1{\myxymessage}

\newdir^{ (}{{}*!/-3pt/\dir^{(}}
\newdir_{ (}{{}*!/-3pt/\dir_{(}}
\newdir^{ )}{{}*!/+3pt/\dir^{)}}
\newdir_{ )}{{}*!/+3pt/\dir_{)}}

\def\rscript#1{\hbox to 0pt{$\scriptstyle#1$\hss}}

\newcommand{\Ab}{\mathop{\CA\it b}\nolimits}

\newcommand{\Cov}{\mathop{Cov^{(\infty)}}\nolimits}
\newcommand{\DM}{\mathop{DM^{eff}_-(k)}\nolimits}
\newcommand{\DeffgM}{\mathop{DM^{eff}_{gm}(k)}\nolimits}
\newcommand{\DgM}{\mathop{DM_{gm}(k)}\nolimits}

\newcommand{\Gnm}{\mathop{\BG_m^n}\nolimits}
\newcommand{\Gnrm}{\mathop{\BG_m^{n-i}}\nolimits}
\newcommand{\Gnrrm}{\mathop{\BG_m^{n-r}}\nolimits}
\newcommand{\imm}{\mathop{{\rm im}}\nolimits}
\newcommand{\ig}{\mathop{i_g}\nolimits}
\newcommand{\isp}{\mathop{i_{\sigma , p}}\nolimits}
\newcommand{\Is}{\mathop{i_\sigma}\nolimits}

\newcommand{\isjZ}{\mathop{i_\sigma^! j_! \, \BZ}\nolimits}

\newcommand{\iTjZ}{\mathop{i_\FT^! j_! \, \BZ}\nolimits}

\newcommand{\iTnmjZ}{\mathop{i_{\FT_{n-m}}^! j_! \, \BZ}\nolimits}

\newcommand{\iYjZ}{\mathop{i_Y^! j_! \, \BZ}\nolimits}
\newcommand{\iYmjZ}{\mathop{i_{Y_m}^! j_! \, \BZ}\nolimits}
\newcommand{\oae}{\mathop{\overline{\{ \eta \}}}\nolimits}
\newcommand{\os}{\mathop{\overline{\{ \sigma \}}}\nolimits}
\newcommand{\oes}{\mathop{\overline{\{ \tau \}}}\nolimits}
\newcommand{\ot}{\mathop{\overline{\{ \nu \}}}\nolimits}
\newcommand{\oit}{\mathop{\overline{\{ \nu_i \}}}\nolimits}

\newcommand{\Schi}{\mathop{Sch^\infty \! /k}\nolimits}
\newcommand{\schi}{\mathop{Sch^{(\infty)} \! /k}\nolimits}
\newcommand{\Spi}{\mathop{S^{\pi (K_1)}}\nolimits}
\newcommand{\Sps}{\mathop{S^{\pi_{[\sigma]} (K_1)}}\nolimits}
\newcommand{\Ses}{\mathop{\FS_{1,[\sigma]}}\nolimits}
\newcommand{\SmC}{\mathop{SmCor(k)}\nolimits}

\newcommand{\Mgm}{\mathop{M_{gm}}\nolimits}
\newcommand{\Mcgm}{\mathop{M_{gm}^c}\nolimits}
\newcommand{\Mcogm}{\mathop{M_{gm}^{\wedge}}\nolimits}
\newcommand{\dMgm}{\mathop{\partial M_{gm}}\nolimits}
\newcommand{\Pes}{\mathop{P_{1,[\sigma]}}\nolimits}
\newcommand{\RC}{\mathop{{\bf R} C}\nolimits}
\newcommand{\ShN}{\mathop{Shv_{Nis}(SmCor(k))}\nolimits}

\newcommand{\uC}{\mathop{\underline{C}}\nolimits}
\newcommand{\uHom}{\mathop{\underline{Hom}}\nolimits}

\newcommand{\Xes}{\mathop{\FX_{1,[\sigma]}}\nolimits}
\newcommand{\Xs}{\mathop{X_\sigma}\nolimits}

\let\oldbullet\bullet
\def\bullet{{\mathchoice{\oldbullet}%
                        {\oldbullet}%
                        {\scriptscriptstyle\oldbullet}%
                        {\oldbullet}}}

\newcommand{\argdot}{{\;\bullet\;}}
\newcommand{\argast}{{\;\ast\;}}

\begin{document}

%

\hfuzz=3pt
\overfullrule=10pt                   


\setlength{\abovedisplayskip}{6.0pt plus 3.0pt}
\setlength{\belowdisplayskip}{6.0pt plus 3.0pt}
\setlength{\abovedisplayshortskip}{6.0pt plus 3.0pt}
\setlength{\belowdisplayshortskip}{6.0pt plus 3.0pt}

\setlength{\baselineskip}{13.0pt}
\setlength{\lineskip}{0.0pt}
\setlength{\lineskiplimit}{0.0pt}

%
%

\title{On the boundary motive of a Shimura variety
\forget{
\footnotemark
\footnotetext{To appear in ...}
}
}
\author{\footnotesize by\\ \\
\mbox{\hskip-2cm
\begin{minipage}{6cm} \begin{center} \begin{tabular}{c}
J\"org Wildeshaus \\[0.2cm]
\footnotesize LAGA\\[-3pt]
\footnotesize Institut Galil\'ee\\[-3pt]
\footnotesize Universit\'e Paris 13\\[-3pt]
\footnotesize Avenue Jean-Baptiste Cl\'ement\\[-3pt]
\footnotesize F-93430 Villetaneuse\\[-3pt]
\footnotesize France\\
{\footnotesize \tt wildesh@math.univ-paris13.fr}
\end{tabular} \end{center} \end{minipage}
\hskip-2cm}
\\[2.5cm]
\forget{
\bf Preliminary version --- not for distribution}\\[1cm]
}
\date{January 22, 2007}
\maketitle
\begin{abstract}
\noindent
Applying the main results of \cite{W9}, we identify the motives
occurring as cones in the co-localization and localization
filtrations of the boun\-dary motive of a Shimura variety. \\

\noindent
Keywords: boundary motives, Shimura varieties, toroidal and 
Baily--Borel compactifications,
completed motives.

\end{abstract}


\bigskip
\bigskip
\bigskip

\noindent {\footnotesize Math.\ Subj.\ Class.\ (2000) numbers:
14G35 (14F42).
}

\eject

\tableofcontents

\bigskip
\vspace*{0.5cm}


%
%

\setcounter{section}{-1}
\section{Introduction}
\label{Intro}



In our previous paper \cite{W9}, we introduced the \emph{boundary
motive} $\dMgm (X)$ of a separated
scheme $X$ of finite type over a perfect field.
There is an exact triangle 
\[
\dMgm (X) \longto \Mgm (X) \longto \Mcgm (X) \longto \dMgm (X) [1] 
\]
in the category of effective motivic complexes, 
relating the boundary motive to $\Mgm (X)$ and $\Mcgm (X)$, 
the \emph{motive} of $X$ and its \emph{motive with compact support}, 
respectively
\cite{VSF}. 
This exact triangle should be employed to construct 
extensions of motives, i.e., classes in motivic cohomology. 
Almost all of the existing attempts to prove the Beilinson or
Bloch--Kato conjectures on special values of $L$-functions
necessitate such classes.
Furthermore, in these cases, the realizations (Betti, de~Rham,
\'etale...) of the classes are constructed
inside the cohomology of (non-compact) Shimura
varieties, using both 
\begin{itemize}
\item[(1)] the respective
realizations of our exact triangle
\[
\dMgm (X) \longto \Mgm (X) \longto \Mcgm (X) \longto \dMgm (X) [1] \; ;
\]
\item[(2)] an accessible description of the realizations of $\dMgm (X)$,
when $X$ is a Shimura variety.
\end{itemize}
This approach is clearly present e.g.\ in Harder's work on 
special values. In particular,
he gave a conjectural description of the \'etale realization of $\dMgm (X)$,
which was eventually proved by Pink \cite{P2}. 
After the work of Harder, Looijenga--Rapoport and Harris--Zucker
\cite{Hd,LR,HrZ},
the Hodge theoretic analogue of this description was later established in \cite{BW}. \\

The aim of this article is to give the motivic version of this description. \\

Let us develop the statement of our main result. 
When the base field admits resolution of singularities, then
there is a \emph{co-localization} principle \cite[Sect.~3]{W9}.
For any compactification $j: X \into \overline{X}$
of $X$, and any 
good stratification $(i_{Y_m}: Y_m \into \overline{X})_m$ of the
boundary $\overline{X} - X$, it 
relates $\dMgm (X)$ to the \emph{motives $\Mgm(Y_m, i_{Y_m}^! j_! \, \BZ)$
of $Y_m$ with
coefficients in $i_{Y_m}^! j_! \, \BZ$},
. When $X$ is smooth, then one also has
\emph{localization} at one's disposal \cite[Sect.~7]{W9},
expressing $\dMgm(X)$ in terms of the
\emph{motives $\Mcgm(Y_m, i_{Y_m}^* j_* \BZ)$ with compact support of $Y_m$ 
and with coefficients in $i_{Y_m}^* j_* \BZ$}. 
Both localization and co-localization should be seen as motivic
analogues of the well-known principles valid for the classical cohomology theories,
relating the cohomology of $\overline{X} - X$ with coefficients in
$j_* \BZ$ to cohomology (with coefficients) of the strata $Y_m$. \\

The main results of this paper identify $\Mgm(Y_m, i_{Y_m}^! j_! \, \BZ)$
(Theorems~\ref{3main} and \ref{4main}) and $\Mcgm(Y_m, i_{Y_m}^* j_* \BZ)$
(Corollaries~\ref{3cor} and \ref{4cor}) when $X$ is a smooth Shimura variety,
and $\overline{X}$ is a toroidal or (when $X$ is pure) the Baily--Borel
compactification, equipped with its canonical stratification. \\

In the toroidal setting, our formulae actually use
little from the material contained in
the present article. We need two main tools:
\emph{analytical invariance} \cite[Sect.~5]{W9} relates the toroidal to a
toric situation. A local computation then identifies the contribution
of a single stratum in a
torus embedding (Lemma~\ref{2A}, Corollary~\ref{2C}). \\

For the purpose of this introduction, let us 
therefore concentrate
on the discussion of the Baily--Borel
compactification $(S^K)^*$ of a pure Shimura variety
$S^K$. Fix a boun\-dary stratum $S_1^K$.
It is (up to an error due to the free action of a finite group) itself
a pure Shimura variety. Denote
by $i$ the immersion of $S_1^K$, and by $j$ the open immersion of $S^K$
into $(S^K)^*$. 
\emph{Invariance under abstract blow-up} \cite[Sect.~4]{W9} shows that one
can compute $M_{gm}(S_1^K, i^! j_! \, \BZ)$ using a
toroidal compactification $S^K (\FS)$ of $S^K$. 
(Invariance under abstract blow-up replaces
the usual application of proper base change in classical cohomology
theories.) Denote the pre-image of $S_1^K$ in
$S^K (\FS)$ by $Z'$. It is then true that the
formal completion of $S^K (\FS)$ along $Z'$ is isomorphic to the
quotient by the free action of a certain arithmetic group $\Delta_1$
of the formal completion of a relative torus embedding along a
union $Z$ of
strata. 
One uses a $\Delta_1$-equivariant version of analytical invariance
(Theorem~\ref{1bK})
in order to identify $M_{gm}(S_1^K, i^! j_! \, \BZ)$
with $R \Gamma (\Delta_1, M_{gm}(Z, i^! j_! \, \BZ)) \,$. Here,
$R \Gamma (\Delta_1,\argdot)$ denotes the ``derived'' functor
of $\Delta_1$-invariants on the level of effective motivic complexes;
we refer to the first half of Section~\ref{1b}
for its definition. 
The problem
is thus reduced to compute the motive with coefficients
in $i^! j_! \, \BZ$ of a union of strata in a torus embedding.
This can be done in some generality (Theorem~\ref{2pC}).
Since the combinatorics of $Z$ is
contractible, the result has a particularly easy shape: 
denoting by $S^{K_1}$ the generic stratum of the relative torus embedding, 
and by $u_1$ the relative dimension of $S^{K_1}$ over the base, one has
\[
M_{gm}(Z, i^! j_! \, \BZ) \cong M_{gm} (S^{K_1}) [u_1]
\]
canonically up to a choice of orientation which will be ignored for
the purpose of this introduction. 
Let us note that $S^{K_1}$
is again a (mixed) Shimura variety. This yields the main isomorphism of Theorem~\ref{4main}:
\bigskip
\[
M_{gm}(S_1^K, i^! j_! \, \BZ) \isoto
R \Gamma (\Delta_1, M_{gm}(S^{K_1})) [u_1]
\]

\bigskip

The results of \cite{BW,P2} give formulae for the degeneration
of mixed (Hodge theoretic or $\ell$-adic) sheaves 
on a Shimura variety. They induce comparison results on the level of
the associated geometric cohomology theories (Betti cohomology in Hodge
theory, \'etale cohomology over the algebraic closure of the reflex field
in the $\ell$-adic context). The mixed sheaves are provided by 
re\-presentations of the group underlying the Shimura variety.
Theorem~\ref{4main} implies
these latter comparison results
for the constant representation $\BZ$ (see Remark~\ref{4rema}). \\
\forget{
Let us add that 
our geometric approach is identical to the one from \cite{BW}
and very near to the one from \cite{P2};
that the main reduction steps of \loccit \ are possible in the motivic
setting results from the tools developed in \cite{W9}. \\
}

Note that at present, the $p$-adic version of the descriptions
from \cite{BW,P2} is not at one's disposal: it would
concern rigid cohomology of closed fibres of a model of the
Baily--Borel compactification
over a ring of integers. However, it is conceivable to deduce
that description from  a suitable version of our main result;
this would necessitate two ingredients, which at present
seem to be within reach: (a)~the construction of a version of 
Voe\-vodsky's motives over Dedekind or at least discrete valuation
rings, (b)~the definition of the rigid realization on this
category. \\

To conclude this introduction,
let us mention that it is possible to ge\-ne\-ralize our main results
to cover the analogues of certain non-constant re\-presentations,
namely those occurring as direct images under the projection from
relative compactifications of semi-Abelian schemes
over the Shimura variety, admitting themselves
an interpretation as Shimura varieties (see Remark~\ref{4rem}).
The additional complication caused by this generalization is
exclusively notational; the technical tools developed in
the present Sections~\ref{1b}
and \ref{2} cover these cases as well. \\

This work was done while I was enjoying a \emph{d\'el\'egation
aupr\`es du CNRS},
and during visits to the \emph{Sonderforschungsbereich~478}
of the University of M\"unster, and to the \emph{Moulin des traits}
at Cinq-Mars-la-Pile. I am grateful to all three institutions.
I also wish to thank D.\ Blotti\`ere and F.\ Lemma for useful discussions
and comments.


\bigskip

%
%

\section{Notations for motives}
\label{1}



The notation of the purely motivic sections of this paper
(Sections~\ref{1a}--\ref{2}) is that
of \cite{W9} (see in particular Section~1 of \loccit), which in turn
follows that of \cite{VSF}. We denote by $k$
the perfect base field over which we work. 
We write $\Schi$ for the category of schemes which are
separated and locally of finite type
over $k$, $Sch/k$ for the full sub-category of objects which are
of finite type
over $k$, and $Sm/k$ for the full sub-category of objects of $Sch/k$ which
are smooth over $k$. Recall the definition
of the category $\SmC$ \cite[p.~190]{VSF}:
its objects are those of $Sm/k$. Morphisms
from $Y$ to $X$ are given by the group $c(Y,X)$ of \emph{finite
correspondences} from $Y$ to $X$.  The category $\ShN$
of \emph{Nisnevich sheaves with transfers}
\cite[Def.~3.1.1]{VSF} is the category of those
contravariant additive functors from $\SmC$ to Abelian groups,
whose restriction to $Sm/k$ is a sheaf for the Nisnevich topo\-logy.
Inside the derived category $D^-(\ShN)$ of complexes bounded from
above, one defines the full triangulated sub-category $\DM$
of \emph{effective motivic complexes} over $k$
\cite[p.~205, Prop.~3.1.13]{VSF} as the one consisting
of objects whose cohomology sheaves are \emph{homotopy invariant}
\cite[Def.~3.1.10]{VSF}.
The inclusion of
$\DM$ into $D^-(\ShN)$ admits a left adjoint $\RC$, which is
induced from the functor
\[
\uC_*: \ShN \longto C^-(\ShN) 
\]
which maps $F$ to the simple complex associated to the 
\emph{singular simplicial complex} \cite[p.~207, Prop.~3.2.3]{VSF}. \\

One defines a functor $L$ from $\Schi$ to
$\ShN$: it associates
to $X$ the Nisnevich sheaf with transfers $c(\argdot,X)$;
note that the above definition of $c(Y,X)$ still makes
sense when $X \in \Schi$ is not
necessarily of finite type and smooth. One defines
the \emph{motive} $\Mgm (X)$ of $X \in \Schi$ as $\RC (L(X))$.
We shall use the same symbol for $\Mgm (X) \in \DM$ and for
its canonical representative $\uC_* (L(X))$ in $C^- (\ShN)$.
There is a second functor $L^c$, defined only for $X \in Sch/k$
as the Nisnevich sheaf of quasi-finite correspondences
\cite[p.~223, 224]{VSF}.
One defines the \emph{motive with compact support}
$\Mcgm (X)$ of $X \in Sch/k$ as
$\RC (L^c(X))$. \\

A second, more geometric approach to motives is
the one developed in \cite[2.1]{VSF}. There, the triangulated
category $\DeffgM$ of \emph{effective geometrical motives} over $k$
is defined. There is a canonical
full triangulated embedding of $\DeffgM$ into $\DM$ \cite[Thm.~3.2.6]{VSF},
which maps the geometrical
motive of $X \in Sm/k$ \cite[Def.~2.1.1]{VSF} to $\Mgm (X)$.
Using this embedding, we consider $\Mgm (X)$ as an object of $\DeffgM$.
Finally, the category $\DgM$ of \emph{geometrical motives} over $k$
is obtained from $\DeffgM$ by inverting the \emph{Tate motive} $\BZ(1)$
\cite[p.~192]{VSF}.
\forget{
All four categories $\DeffgM$, $\DgM$, $D^-(\ShN)$, and $\DM$
are tensor triangulated, and admit unit objects, which we denote
by the same symbol $\BZ (0)$
\cite[Prop.~2.1.3, Cor.~2.1.5, p.~206,
Thm.~3.2.6]{VSF}.
For $M \in \DgM$ and $n \in \BZ$, write $M(n)$ for the tensor product
$M \otimes \BZ (n)$.
}
According to \cite[Thm.~4.3.1]{VSF},
the functor $\DeffgM \to \DgM$ is a full triangulated
embedding if $k$ admits resolution of singularities.


\bigskip

%
%

\section{The completed motive}
\label{1a}



The torus embeddings we need to consider later
will in general only be \emph{locally} of finite type. We are thus led to define a new
extension of the functor 
\[
\Mgm: Sch/k \longto \DM
\]
to a certain sub-category of $\Schi$.
We introduce the notion of \emph{completed motive}.
It will turn out
to be very well behaved with respect to free actions of possibly
infinite abstract groups, to be considered in the next section. 
Also, its definition will allow for a straight-forward
modification of analytical invariance to the geometric context
we need to consider later. 

\begin{Def} \label{1aA}
Let $X \in \Schi$. \\[0.1cm]
(a)~Denote by $\Cov (X)$ the system of Zariski coverings 
$\FU = (U_\alpha)_\alpha$
of $X$ such that all $U_\alpha$ are in $Sch/k$, and such that for 
all open subsets
$V \subset X$ belonging to $Sch/k$, the set
$\{ \alpha \tei U_\alpha \cap V \ne \emptyset \}$
is finite. \\[0.1cm]
(b)~Denote by $\schi$ the full sub-category of $\Schi$ of objects $X$ for which
$\Cov (X)$ is non-empty.
\end{Def}

To the knowledge of the author, the class $\schi$ of schemes has not been
considered previously.

\begin{Ex} \label{1aEx}
Examples of schemes in $\Schi$ not belonging to $\schi$ occur in the context of
torus embeddings (see Section~\ref{2}). Let $T$ be a split torus over $k$, and 
$\FS$ a rational partial polyhedral
decomposition (see e.g.\ \cite[Chap.~I]{KKMS}, or \cite[Chap.~5]{P1})
of $Y_*(T)_\BR$, where $Y_*(T)$ denotes the cocharacter group of $T$.
It gives rise to an open dense embedding of $T$ into a $k$-scheme $T_\FS$ 
which is locally of finite type. Since $T$ meets any open subset of 
the torus embedding $T_\FS$, this scheme lies in $\schi$ if and only if
$\FS$ is finite.   
\end{Ex}

Note that for $X \in \schi$, one can find $\FU \in \Cov (X)$ such that
all members of $\FU$ are affine.

\begin{Def} \label{1aF}
Let $Y \into W$ be a closed immersion in $\Schi$. \\[0.1cm]
(a)~Define
\[
L(W / W-Y) := L(W) / L(W-Y) \; ,
\]
and
\[
\Mgm (W / W-Y) := \RC (L (W / W - Y)) \; .
\]
(b)~Let $\FU_W$ be a Zariski covering of an open subset $V$ of
$W$, and $X$ a sub-scheme of $W$.
Denote by $\FU_X$ the induced covering of $V \cap X$. \\[0.1cm]
(c)~Denote by $\Cov (W)_Y$ the system of 
Zariski coverings $\FU_W$ of open subsets of
$W$ such that $\FU_Y$ is in $\Cov (Y)$. (In particular, $\FU_W$
covers the whole of $Y$.)
\end{Def}

Note that by definition, $\Cov (W)_Y$ is
non-empty if and only if $Y$ belongs to $\schi$.

\begin{Rem} \label{1aG}
Let $Y \into W$ be a closed immersion in $\Schi$, and 
$\FU_W = (U_\alpha)_\alpha$ a Zariski covering of an open subset
of $W$ containing $Y$. Then
$\FU_W$ belongs to $\Cov (W)_Y$ if and only if for all $\beta$, the set
$\{ \alpha \tei U_\alpha \cap U_\beta \not \subset W - Y \}$
is finite.
\end{Rem}

For $Y \into W$ as in Definition~\ref{1aF} and $\FU_W = (U_\alpha)_\alpha \in \Cov (W)_Y$,
construct the \v{C}ech complex $L^\bullet (\FU_W / \FU_{W - Y})$,
concentrated in non-positive degrees, as follows: its $-m$-th term
is the direct product of the $L(U_A / U_A -Y)$, 
for all $(m+1)$-fold intersections $U_A$
of the $U_\alpha$. Note that the differentials are 
defined thanks to Remark~\ref{1aG}.

\begin{Lemma} \label{1aH}
Let $Y \into W$ be a closed immersion in $\Schi$ and
$\FU_W  = (U_\alpha)_{\alpha \in I} \in \Cov (W)_Y$. \\[0.1cm]
(a)~If $\FU_Y$ is finite (and hence $Y$ is in $Sch/k$), then the augmentation
\[
L^\bullet (\FU_W / \FU_{W - Y}) \longto L(W / W-Y)
\]
is a quasi-isomorphism in $C^- (\ShN)$. \\[0.1cm]
(b)~Assume that $\FV_W  = (V_\beta)_{\beta \in J} \in \Cov (W)_Y$ 
is a refinement of $\FU_W$:
\[
\forall \, \beta \in J \, \exists \, \alpha \in I \, ,
\; V_\beta \subset U_\alpha \; .
\] 
Then any refinement map $\iota: J \to I$ induces a quasi-isomorphism
\[
L^\bullet (\iota): L^\bullet (\FV_W / \FV_{W - Y}) \longto 
L^\bullet(\FU_W / \FU_{W - Y})
\]
in $C^- (\ShN)$. Different choices of refinement maps $\iota$ lead to
homotopy equivalent morphisms $L^\bullet (\iota)$.
\end{Lemma}

\begin{Proof}
The very last statement is of course well known.
Note that for
$W' \in \Schi$, the sheaf $L(W')$ is the filtered direct limit of the
$L(X)$, with $X$ running through the open sub-schemes of $W'$
which are of finite type over $k$. Hence \cite[Prop.~3.1.3]{VSF}
continues to hold for Nisnevich coverings of schemes in $\Schi$.
Part~(a) follows from this --- observe that since $\FU_Y$ is finite,
the direct product $\prod_{\alpha \in A} L(U_\alpha / U_\alpha - Y)$ 
equals $L( \coprod_{\alpha \in A} U_\alpha ) / 
L( \coprod_{\alpha \in A} U_\alpha - Y)$, where $A$ denotes
the set of indices $\alpha$ for which $U_\alpha \not \subset W-Y$.
In order to prove (b), define a third covering
\[
\FW_W := (U_\alpha \cap V_\beta)_{\alpha, \beta} \; .
\]
It is easy to see that it also belongs to $\Cov (W)_Y$. Furthermore, there are
canonical refinement maps from $L^\bullet(\FW_W / \FW_{W-Y})$ to 
$L^\bullet(\FU_W / \FU_{W-Y})$ and to
$L^\bullet(\FV_W / \FV_{W-Y})$. That these
are quasi-isomorphisms follows from
(a), applied to each component of $L^\bullet(\FU_W / \FU_{W-Y})$ and 
$L^\bullet(\FV_W / \FV_{W-Y})$,
respectively. Finally, observe that the diagram
\[
\xymatrix@R-20pt{
& L^\bullet(\FW_W / \FW_{W-Y}) \ar[ddl] \ar[ddr] & \\
& & \\
L^\bullet(\FV_W / \FV_{W-Y}) \ar[rr] & & L^\bullet(\FU_W / \FU_{W-Y})
\\}
\]
is commutative up to homotopy.
\end{Proof}

Hence the following definition makes sense:

\begin{Def} \label{1aI}
Let $Y \into W$ be a closed immersion in $\Schi$, with $Y \in \schi$.
The \emph{completed motive of $W$ relative to $W - Y$} is defined as
\[
\Mcogm (W / W - Y) := 
\invlim_{\FU_W \in \Cov (W)_Y} \RC (L^\bullet (\FU_W / \FU_{W - Y})) \; .
\]
\end{Def}

Given that all transition maps in the inverse limit are isomorphisms,
our definition should be seen as an attempt to avoid making choices.
Of course, $\Mcogm (W / W - Y)$ is represented by the
simple complex associated to $\uC_* (L^\bullet (\FU_W / \FU_{W-Y}))$, for any
$\FU_W \in \Cov (W)_Y$. In particular:

\begin{Prop} \label{1aJ}
Let $Y \into W$ be a closed immersion in $\Schi$, with $Y \in Sch/k$. 
Then the augmentation
\[
\Mcogm (W / W - Y) \longto \Mgm (W / W - Y)
\]
is an isomorphism in $\DM$.
\end{Prop}

The case $Y = W$ is important; let us introduce a specific notion:

\begin{Def} \label{1aD}
Let $Y \in \schi$. The \emph{completed motive of $Y$} is defined as
\[
\Mcogm (Y) := \Mcogm (Y / \emptyset) \; .
\]
\end{Def}

We want to 
consider the geometric situation in which \emph{co-localization} works
\cite[Section~3]{W9}: fix closed immersions $Y \into Y' \into W$ in $\Schi$.
Write $j$ for the open immersion of $W-Y'$, and $i_Y$ for the closed
immersion of $Y$ into $W$. Denote by $\FY$ the commutative diagram
\[
\xymatrix@R-20pt{
Y' - Y \ar[r] \ar[dd] &
W - Y \ar[dd] \\
& \\
Y' \ar[r] &
W
\\}
\]
Assume in addition that $Y \in \schi$. For $\FU_W \in \Cov (W)_Y$, 
denote by $L^\bullet (\FU_{\FY})$ the simple complex
associated to
\[
\xymatrix@R-20pt{
L^\bullet (\FU_{Y' - Y}) \ar[r] \ar[dd] &
L^\bullet (\FU_{W - Y}) \ar[dd] \\
& \\
L^\bullet (\FU_{Y'}) \ar[r] &
L^\bullet (\FU_W)
\\}
\]
with $L^\bullet (\FU_W)$ sitting in degree $(0,\argast)$.
As before, one shows that the following is well-defined:

\begin{Def} \label{1aK}
The \emph{completed motive of $Y$ with
coefficients in $\iYjZ$} is defined as
\[
\Mcogm (Y, \iYjZ) :=
\invlim_{\FU_W \in \Cov (W)_Y} \RC (L^\bullet (\FU_{\FY})) \; .
\]
\end{Def}

The completed motive $\Mcogm (Y, \iYjZ)$
is represented by the
simple complex associated to $\uC_* (L^\bullet (\FU_{\FY}))$, for any
$\FU_W \in \Cov (W)_Y$.

\begin{Prop} \label{1aL}
In the situation of the definition, there is a canonical exact triangle
\[
\Mcogm (Y, \iYjZ)[-1] \to \Mcogm (Y' / Y' - Y) \to
\Mcogm (W / W - Y) \to \Mcogm (Y, \iYjZ)
\]
in $\DM$. If $Y \in Sch/k$, then the natural morphism
\[
\Mcogm (Y, \iYjZ) \longto \Mgm (Y, \iYjZ)
\]
is an isomorphism in $\DM$. Here, the right hand side denotes the
motive of $Y$ with coefficients in $\iYjZ$ \cite[Def.~3.1]{W9}.
\end{Prop}

\begin{Rem}
Definition~\ref{1aK} and \cite[Def.~3.1]{W9}
extend to situations where the immersion
of $Y$ into $W$ (hence into $Y'$) is only locally closed:
one replaces $W$ by an open subset containing $Y$ as a closed
sub-scheme.
Proposition~\ref{1aL} continues to hold in this larger
generality. We leave the details to the reader.
\end{Rem}

\bigskip

%
%

\section{Equivariant completed motives}
\label{1b}



For the application to Shimura varieties, we need to consider the situation
of the preceding section under the additional hypothesis of the presence
of the action of an abstract group $H$. 
The result we are aiming at is an $H$-equivariant version of 
\emph{analytical invariance} \cite[Thm.~5.1]{W9}.
Let us start by recalling some basic material on the right derived functor
of the fixed point functor associated to $H$. \\

Fix an Abelian category $\CA$ and an abstract group $H$.
We assume that (A)~$\CA$ is closed under arbitrary products,
and that (B)~the trivial $\BZ H$-module $\BZ$ admits a bounded resolution
\[
0 \longto F_{n} \longto \ldots \longto
  F_{1}\longto F_{0}\longto \BZ \longto 0
\]
by free $\BZ H$-modules.

\begin{Rem} \label{1bA}
Recall that $H$ is said to be of type $FL$ if in addition the $F_i$ in
the above resolution can be chosen to be of finite type over $\BZ H$.
\end{Rem}

We shall denote by $H \text{-}\CA$ the Abelian category
of objects of $\CA$ provided with a left $H$-action.
If $\gamma \in H$ and
$A\in H \text{-}\CA$, we denote by the same letter
$\gamma$ the corresponding automorphism of $A$. 
Morphisms in $H \text{-}\CA$ between two objects $A$ and $B$ 
are those morphisms $f: A \to B$ in $\CA$ commuting with the action
of $H$, i.e., such that 
\[
f \circ \gamma = \gamma \circ f : A \longto B \; , \; \forall \; 
\gamma \in H \; .
\]
We denote by $e$
the unit element of $H$.
By definition,  the fixed point functor associated to $H$
is the functor
\[
\Gamma (H,\argdot) = (\argdot)^{H}:
H \text{-}\CA \longrightarrow \CA
\]
given by
\[
\Gamma (H,A) = (A)^{H} := \bigcap_{\gamma \in H}\kernel(e-\gamma) \; .
\]
The condition (A) we put on $\CA$ implies that
the right derived functor
\[
R\Gamma (H,\argdot) : D^{+}(H \text{-}\CA)
\longto D^{+}(\CA)
\]
exists (e.g. \cite[Thm.~3.11, Var.~3.18~(a)]{BW}). 
We need to be more precise. Denote by $\Ab$ the category of Abelian groups.
Recall the definition of the functor
\[
\Hom: \Ab \times \CA \longto \CA 
\]
\cite[Sect.~3]{BW}. For $(M,A) \in \Ab \times \CA$, the object
$\Hom (M,A)$ represents the functor associating to $B$ all group
morphisms from $M$ to $\Hom_{\CA} (B,A)$. 
If $M$ is free with basis $(x_i)_{i \in I}$, then
\[
\Hom (M,A) \cong \prod_{i \in I} A \; .
\]
This construction admits an extension to complexes; we denote it
by the same symbol
\[
\Hom: C^b (\Ab) \times C^?(\CA) \longto C^?(\CA) \; ,
\]
for $? \in \{ b, +, -, \text{blank} \}$.
Homotopy equivalent complexes $M^1_\ast$ and $M^2_\ast$ in
$C^b (\Ab)$ yield homotopy equivalent complexes 
\[
\Hom (M^1_\ast,A^\bullet) \quad \text{and} \quad
\Hom (M^2_\ast,A^\bullet)
\]
in $C(\CA)$, for any $A^\bullet \in C(\CA)$. 
The same statement is true for the variant
\[
\Hom: C^+ (\Ab) \times C^-(\CA) \longto C^-(\CA) 
\]
of the functor $\Hom$. \\
 
For $(M_\ast,A^\bullet) \in C^b (H \text{-} \Ab) \times C(H \text{-}\CA)$,
there is the diagonal action of $H$ on
the complex $\Hom (M_\ast,A^\bullet)$ \cite[Def.~3.12]{BW}.  
The sub-complex of (componentwise) invariants
is denoted by $\Hom (M_\ast,A^\bullet)^H$;
note that this complex is in $C^? (\Ab)$ if $A^\bullet$ is
in $C^? (H \text{-} \Ab)$.
We then have \cite[Prop.~3.13~(b), Var.~3.18~(a)]{BW}:

\begin{Prop} \label{1bB}
Let $F_\ast \to \BZ$ be any bounded resolution of
the trivial $\BZ H$-module $\BZ$ by free
$\BZ H$-modules. Then the functor 
\[
R\Gamma (H,\argdot) : D^{+}(H \text{-}\CA)
\longto D^{+}(\CA)
\]
is induced by $\Hom(F_{\ast},\argdot)^{H}$.
\end{Prop}
 
This conclusion would still hold without the boundedness assumption
on $F_\ast$. However,
we are interested in the following variant, for which the condition (B) 
we put on $H$ is essential:

\begin{Cor} \label{1bC}
Let $F_\ast \to \BZ$ be any bounded resolution of
the trivial $\BZ H$-module $\BZ$ by free
$\BZ H$-modules. Then the right derived functors 
\[
R\Gamma (H,\argdot): D^b(H \text{-}\CA)
\longto D^b(\CA) \; ,
\]
\[
R\Gamma (H,\argdot): D(H \text{-}\CA)
\longto D(\CA) 
\]
and 
\[
R\Gamma (H,\argdot): D^-(H \text{-}\CA)
\longto D^-(\CA) 
\]
of $\Gamma (H,\argdot)$ exist and are induced by 
the functors
\[
\Hom(F_{\ast},\argdot)^{H}: C^?(H \text{-}\CA) 
\longto C^?(\CA) \; .
\]
In particular, all the functors denoted $R\Gamma (H,\argdot)$
are compatible under restriction.
\end{Cor}

\begin{Proof}
Proposition~\ref{1bB} implies that the functor $\Gamma (H,\argdot)$
is of bounded cohomological dimension.
Our claim thus follows from \cite[Thm.~II.2.2.2, Cor.~2]{V}. 
\end{Proof}

\begin{Cor} \label{1bD}
(a)~Let $? \in \{ b, +, -, \text{blank} \}$. The functor
\[
R\Gamma (H,\argdot): D^?(H \text{-}\CA) 
\longto D^?(\CA)
\]
is right adjoint to the functor $\Res^1_H$ associating to a complex of objects
of $\CA$ the same complex with trivial $H$-action. \\[0.1cm]
(b)~Let $M$ be a free $\BZ H$-module, and $A \in H \text{-}\CA$.
Then $\Hom(M,A)$ is  
$\Gamma (H,\argdot)$-acyclic.
\end{Cor}

\begin{Proof}
Part~(a) results from the fact that the corresponding statement is
true for the functor $\Gamma (H,\argdot)$ on $C^?(H \text{-}\CA)$,
and from the description of $R\Gamma (H,\argdot)$ from Corollary~\ref{1bC}. 
For (b), choose $F_\ast \to \BZ$ as in \ref{1bC}.
Observe that by adjunction,
\[
\Hom(F_{\ast},\Hom(M,A))^H = 
\Hom(F_{\ast} \otimes_\BZ M,A)^H \; .
\]
The complex $F_{\ast} \otimes_\BZ M$ is a free resolution of the free
$\BZ H$-module $M$. Therefore, the right hand side is homotopy
equivalent to $\Hom(M,A)^H$, a complex concentrated in degree zero.
\end{Proof} 

\begin{Rem} \label{1bDa}
Applying the dual of the above method, one shows the existence of the
left derived functor
\[
L \Lambda (H,\argdot): D^-(H \text{-}\CB)
\longto D^-(\CB) 
\]
of the functor of co-invariants
\[
\Lambda (H,\argdot) = (\argdot)_{H}:
H \text{-}\CB \longrightarrow \CB
\]
for any Abelian category $\CB$ closed under arbitrary sums. 
It is left adjoint to the functor $\Res^1_H$. Note that
we do not need to impose condition (B) on the group $H$ as long as we
content ourselves to work with complexes bounded from above. 
\end{Rem}

Now we specialize to the case $\CA = \ShN$. We have
\[
R\Gamma (H,\argdot): D^-(H \text{-}\ShN) 
\longto D^-(\ShN) \; ,
\]
and this functor respects the full sub-categories $D^b$.

\begin{Def} \label{1bE}
(a)~An object $F$ of $H \text{-}\ShN$ is called \emph{homotopy invariant}
if the underlying object of $\ShN$ is homotopy invariant. \\[0.1cm]
(b)~The category $H \text{-} \DM$ 
of \emph{effective $H$-equivariant motivic complexes} over $k$
is the full sub-category of
$D^-(H \text{-}\ShN)$ consisting of objects whose cohomology sheaves
are homotopy invariant.
\end{Def}

One shows \cite[Cor.~3.5]{S} that one can argue as in the non-equivariant
situation \cite[3.1]{VSF}, to see that $H \text{-} \DM$ is a
triangulated sub-category of $D^-(H \text{-}\ShN)$.
It follows from Corollary~\ref{1bC} that $R \Gamma (H, \argdot)$ induces
\[
R \Gamma (H, \argdot): H \text{-} \DM \longto \DM \; .
\]
The functor 
\[
\uC_*: \ShN \longto C^-(\ShN) 
\]
admits an obvious $H$-equivariant version, and induces
\[
\RC: D^-(H \text{-}\ShN) \longto H \text{-} \DM \; .
\]
One shows \cite[Thm.~3.7]{S} that as in the non-equivariant situation,
$\RC$ is left adjoint to the
inclusion of $H \text{-} \DM$ into $D^-(H \text{-}\ShN)$.
Furthermore:

\begin{Prop} \label{1bF}
The diagram 
\[
\xymatrix@R-20pt{
D^-(H \text{-}\ShN) \ar[rr]^-{R \Gamma (H, \argdot)} \ar[dd]_{\RC} & &
D^-(\ShN) \ar[dd]^{\RC} \\
 & & \\
H \text{-} \DM \ar[rr]^-{R \Gamma (H, \argdot)} & &
\DM
\\}
\]
commutes.
\end{Prop}

\begin{Proof}
By definition, the functors $\uC_*$ and
$\Hom (F_\ast,\argdot)^H$ commute.
\end{Proof}

\begin{Rem} \label{1bFa}
Dually, we get (even without condition (B) on $H$)
the functor 
\[
L \Lambda (H,\argdot): D^-(H \text{-}\ShN)
\longto D^-(\ShN) 
\]
restricting to  
\[
L \Lambda (H,\argdot): H \text{-} \DM \longto \DM
\]
and commuting with $\RC$. 
\end{Rem}

Denote by $H \text{-} \Schi$ the category
of objects of $\Schi$ provided with a left $H$-action; similarly for
$H \text{-} \schi$. Observe that for $X \in H \text{-} \Schi$, the sheaf $L(X)$
is canonically endowed with an $H$-action.

\begin{Def} \label{1bG}
Let $X \in H \text{-} \Schi$. Define the \emph{$H$-equivariant motive}
$\Mgm (X) \in H \text{-} \DM$ as the image of $L(X) \in H \text{-}\ShN$
under $\RC$.
\end{Def}

Of course, the natural forgetful functor maps the $H$-equivariant
motive to the usual motive. 

\begin{Rem} \label{1bGa}
One can (but we shall not in this article) 
consider the full triangulated sub-category
$H \text{-} \DeffgM$ of $H \text{-} \DM$ of
\emph{effective $H$-equivariant geometrical motives} over $k$,
i.e., the full triangulated sub-category generated by the $\Mgm (X)$,
for $X \in H \text{-} Sm/k$.  
The category of \emph{$H$-equivariant geometrical motives} 
$H \text{-} \DgM$ is obtained from $H \text{-} \DeffgM$
by inverting $\Res^1_H \BZ (1)$.
If $H$ is of type $FL$, then both $R\Gamma (H,\argdot)$ and 
$L \Lambda (H,\argdot)$ send $H \text{-} \DeffgM$
to $\DeffgM$. Furthermore, by the universal property of 
$H \text{-} \DgM$, they induce functors
\[
R\Gamma (H,\argdot) \; , \; L \Lambda (H,\argdot):
H \text{-} \DgM \longto \DgM 
\]
right resp.\ left adjoint to $\Res^1_H$. Assume that
the category $H \text{-} \DgM$ were known to be a rigid tensor
category, in a way compatible with the tensor structure on $\DgM$
under $\Res^1_H$. Then duality in $H \text{-} \DgM$ and $\DgM$
would exchange $R\Gamma (H,\argdot)$ and $L \Lambda (H,\argdot)$
because they are the adjoints of the same tensor functor.
\end{Rem}

\begin{Def} \label{1bH}
Let $X \in H \text{-} \Schi$, and assume that the action of $H$ 
on $X$ is faithful. \\[0.1cm]
(a)~The action of $H$ on $X$ is \emph{proper} if there is a Zariski
covering $\FU$ of $X$ by open affines
such that for any $U$ in $\FU$, the set
$\{ \gamma \in H \tei U \cap \gamma U \ne \emptyset \}$ is finite. \\[0.1cm]
(b)~The action of $H$ on $X$ is \emph{free} if the stabilizer of any
point of $X$ is reduced to $\{ e \}$.
\end{Def}

Note that if $H$ acts properly on $X \in \Schi$, 
then the quotient $H \backslash X$ can be formed. 
Furthermore,
if $U_1, U_2 \in Sch/k$ are open sub-schemes of $X$, then the set
$\{ \gamma \in H \tei U_1 \cap \gamma U_2 \ne \emptyset \}$ is finite. 
If in addition the
action of $H$ is free,
then there is a Zariski
covering $\FU$ of $X$ by open affines
such that for any $U$ in $\FU$, we have
\[
\{ \gamma \in H \tei U \cap \gamma U \ne \emptyset \} = \{ e \} \; .
\]
Hence in this case, the projection $\Pi: X \to H \backslash X$
is a local isomorphism, and $\FU$ induces a Zariski covering 
$\Pi (\FU)$ of $H \backslash X$ by open affines. 
By definition, a covering $\FV$ of $H \backslash X$ by affines splits $\Pi$
if there is $\FU$ as above such that $\FV = \Pi (\FU)$.
We leave the proof of the following to the reader:

\begin{Prop} \label{1bI}
Let $X \in \Schi$, with a free and proper action of $H$. 
Assume that $H \backslash X \in \schi$. \\[0.1cm]
(a)~There exists a covering $\FV \in \Cov (H \backslash X)$ by affines
splitting $\Pi$. \\[0.1cm]
(b)~$X \in \schi$.
\end{Prop} 

We are ready to state the main result of this section. 
Fix closed immersions $Y_1 \into W_1$ in $H \text{-}\Schi$ and
$Y_2 \into W_2$ in $\Schi$, such that $Y_2 \in \schi$, and such that
the action of $H$ on $Y_1$ is free and proper.
Assume given an isomorphism $f: Y_2 \isoto H \backslash Y_1$,
which extends to
an isomorphism
\[
f: (W_2)_{Y_2} \isoto H \backslash (W_1)_{Y_1} \; ,
\]
where $(W_m)_{Y_m}$ denotes the 
formal completion of $W_m$ along $Y_m$, $m = 1,2$.
Note that since the action of $H$ on $Y_1$ is free and proper,
the same is true for the action on $(W_1)_{Y_1}$. In particular,
the quotient $\Pi: (W_1)_{Y_1} \to H \backslash (W_1)_{Y_1}$ can be
formed and is a local isomorphism. Furthermore, Proposition~\ref{1bI}~(b)
implies that $Y_1 \in \schi$. Using coverings of $Y_1$ induced by $f$ and
elements of $\Cov (Y_2)$, one sees that 
$\Mcogm (W_1 / W_1 - Y_1) \in \DM$ has a canonical pre-image
under the forgetful functor from $H \text{-}\DM$. We use the same
symbol $\Mcogm (W_1 / W_1 - Y_1)$ for this pre-image.

\begin{Thm} \label{1bmain}
Pull-back of finite correspondences induces an isomorphism
\[
\Mcogm (W_2 / W_2 - Y_2) \isoto 
R\Gamma(H, \Mcogm (W_1 / W_1 - Y_1)) 
\]
in $\DM$, and depending only on $f$. 
\end{Thm}

\begin{Proof}
Use $f$ to identify $Y_2$ and $H \backslash Y_1$. Proposition~\ref{1bI}~(a)
enables us to choose a covering $\FU_2 \in \Cov (Y_2)$ by affines such that
$\Pi: Y_1 \to Y_2$ is split over $\FU_2$. Define $\FU_1 \in \Cov (Y_1)$
as the affine covering induced by $\FU_2 \,$: for any $U_2 \in \FU_2$, its
pre-image $\Pi^{-1} (U_2)$ is thus the disjoint union of $|H|$ affines,
each of which is isomorphic to $U_2$ and belongs to $\FU_1$. 
By definition of the completed motive, $\Mcogm (W_m / W_m - Y_m)$
is represented by 
\[
\uC_* (L^\bullet (\FU_{W_m} / \FU_{W_m - Y_m})) \; ,
\]
for any Zariski covering $\FU_{W_m}$ of an open subset of $W_m$ such that
$\FU_{Y_m} = \FU_m$, $m = 1,2$. Using \cite[Thm.~5.5]{W9} (including the
statements (b) and (c) on compatibility under restriction), one sees that
pull-back under $\Pi$ induces canonically
\[
L^\bullet (\FU_{W_2} / \FU_{W_2 - Y_2}) \isoto
L^\bullet (\FU_{W_1} / \FU_{W_1 - Y_1})^H \; .
\]
\emph{Choose} for each $U_2$ in $\FU_2$ a $U_1$ in $\FU_1$
contained in its pre-image under $\Pi$. This choice induces,
again by \cite[Thm.~5.5]{W9},
an isomorphism 
\[
\Hom(\BZ H, L^i (\FU_{W_2} / \FU_{W_2 - Y_2})) 
\isoto L^i (\FU_{W_1} / \FU_{W_1 - Y_1}) \; ,
\]
for each $i \in \BZ$.
Here, $L^i (\FU_{W_2} / \FU_{W_2 - Y_2})$ carries
the trivial action of $H$.
By Corollary~\ref{1bD}~(b), the components of 
$L^\bullet (\FU_{W_1} / \FU_{W_1 - Y_1})$ are therefore  
$\Gamma (H,\argdot)$-acyclic. Hence its $H$-invariants represent
\[
R\Gamma(H, L^\bullet (\FU_{W_1} / \FU_{W_1 - Y_1})) \; .
\]
This shows that we even have an isomorphism
\[
L^\bullet (\FU_{W_2} / \FU_{W_2 - Y_2}) \isoto
R\Gamma(H, L^\bullet (\FU_{W_1} / \FU_{W_1 - Y_1})) \quad \quad \quad (\ast)
\]
in $D^- (\ShN)$. Our claim then follows from Proposition~\ref{1bF}.
\end{Proof}

We add two geometrical data:
immersions $Y_1 \into Y'_1 \into W_1$ in $H \text{-}\Schi$ and
$Y_2 \into Y'_2 \into W_2$ in $\Schi$, with $Y_2 \in \schi$, such that
$Y'_m \into W_m$ are closed, and 
the action of $H$ on $Y_1$ is free and proper. Write
$j_m$ for the open immersion of $W_m-Y'_m$, and $i_{Y_m}$ for the 
immersion of $Y_m$ into $W_m$. 
As above, one notes that there is a canonical pre-image
$\Mcogm (Y_1, i_{Y_1}^! j_{1 !} \, \BZ)$ in $H \text{-}\DM$
of the completed motive of $Y_1$ with
coefficients in $i_{Y_1}^! j_{1 !} \, \BZ$. 
We have the following $H$-equivariant version of \cite[Thm.~5.1]{W9}:

\begin{Thm} \label{1bK}
Assume gi\-ven an isomorphism $f: Y_2 \isoto H \backslash Y_1$,
which extends to an isomorphism
\[
f: (W_2)_{Y_2} \isoto H \backslash (W_1)_{Y_1} 
\]
inducing an isomorphism $(Y'_2)_{Y_2} \cong H \backslash (Y'_1)_{Y_1}$.
Then $f$ induces an isomorphism 
\[
\Mcogm (Y_2, i_{Y_2}^! j_{2 !} \, \BZ) \isoto 
R\Gamma ( H,\Mcogm (Y_1, i_{Y_1}^! j_{1 !} \, \BZ)) 
\]
in $\DM$.
\end{Thm}

\begin{Proof}
This follows from the fact that the isomorphism $(\ast)$ in
the proof of Theorem~\ref{1bmain} is compatible with restriction
of the $W_m$ to $Y'_m$ \cite[Thm.~5.5~(b)]{W9}.
\end{Proof}


\bigskip

%
%

\section{Motives of torus embeddings}
\label{2}



In this section, we shall consider the following situation:
$k$ is a perfect base field, $T$ a split 
torus over $k$, $B \in Sch/k$ a reduced scheme,
and $X \to B$ a $T$-torsor.
We fix a \emph{rational partial polyhedral
decomposition} (see e.g.\ \cite[Chap.~I]{KKMS}, or \cite[Chap.~5]{P1})
$\FS$ of $Y_*(T)_\BR$,
where $Y_*(T)$ denotes the cocharacter group of $T$, and a
non-empty proper subset $\FT \subset \FS$ satisfying the following condition:
$(\alpha)$~whenever $\tau \in \FT$ is a face of $\sigma \in \FS$, then 
$\sigma$ belongs
to $\FT$. Recall that by definition, any element of $\FS$ is a 
\emph{convex rational
polyhedral cone} in $Y_*(T)_\BR$, any face of any cone in $\FS$, and
the intersection of any two cones in $\FS$ belong to $\FS$, and 
$\{ 0 \} \in \FS$. This implies that $0 := \{ 0 \}$ is a face of any cone in
$\FS$, and hence that no such cone contains a non-trivial linear subspace
of $Y_*(T)_\BR$. Since $\FT$ is supposed to be a proper subset of $\FS$,
the cone $0$ does not belong to $\FT$. We assume that $(\beta)$~any cone
$\sigma$ of $\FS$ is \emph{smooth}, i.e., the semi-group $\sigma \cap Y_*(T)$
can be generated by a subset of a $\BZ$-basis of $Y_*(T)$. \\

To the given data are associated, as is described e.g.\ in \cite[5.5]{P1},
(1)~a relative \emph{torus embedding} $j:X \into X_\FS$ 
identifying $X$ with an open dense sub-scheme of a scheme $X_\FS$,
which is locally of finite type over $B$,
(2)~a stra\-ti\-fication of $X_\FS$ into locally closed 
reduced sub-schemes $\Is: \Xs \into X_\FS$
indexed by the cones $\sigma \in \FS$,
and of finite type over $B$.
This stratification allows to 
associate to any \emph{locally
closed} subset $\FV$ of $\FS$ the locally closed reduced sub-scheme
$X_\FV$ of $X_\FS$ which set-theoretically is the union  
of the $\Xs$, with $\sigma \in \FV$; here, a subset $\FV$ of $\FS$
is called locally closed if the relations $(\nu \in \FV$ is a face of
$\sigma \in \FS)$ and $(\sigma$ is a face of $\mu \in \FV)$ imply 
that $\sigma \in \FV$. 
The scheme $X_\FV$ is of finite type over $B$
if and only if $\FV$ is finite.
In particular, the above date define
(3)~a sub-scheme $i_\FT: X_\FT \into X_\FS$ contained in the complement
of $X = X_0$. Note that $X_\FT$ 
is closed in $X_\FS$
because of condition~$(\alpha)$, and that $X_\FS$ is smooth over $B$
because of condition~$(\beta)$. \\

We add one more condition on the given data. In order to formulate
it, define
\[
D := \bigcup_{\tau \in \FT} \tau^\circ \subset Y_*(T)_\BR \; ,
\]
where for each cone $\tau$ we denote by $\tau^\circ$ the topological
interior of $\tau$ inside the linear subspace of $Y_*(T)_\BR$ generated by
$\tau$. The subset $D$ of $Y_*(T)_\BR$ is endowed with the induced
topology. Note that $0 \not \in D$ since $0 \not \in \FT$, and
that $D$ is \emph{conic} in the sense that multiplication
by strictly positive real numbers on $Y_*(T)_\BR$ preserves $D$. 
We assume that $(\gamma)$~every point of $D$ admits
a neighbourhood $U$ such that $U \cap \tau \ne \emptyset$ for only a
finite number of $\tau \in \FT$. 
Note that $(\gamma)$ implies that
every $\tau \in \FT$ is a face of only finitely
many $\sigma \in \FT$. Hence the closure of any 
$X_\tau$, for $\tau \in \FT$,
is of finite type over $B$. Furthermore, writing
\[
\os := \{ \nu \in \FS \tei \nu \; \text{face of} \; \sigma \}
\]
for any cone $\sigma$ of $\FS$, we have:

\begin{Prop} \label{2pnull}
(a)~The scheme $X_\FT$ lies in $\schi$. \\[0.1cm]
(b)~The system $(X_{\oes})_{\tau \in \FT}$ belongs to $\Cov (X_\FS)_{X_\FT}$.
\end{Prop}

This result should be compared to Example~\ref{1aEx}. \\

\begin{Proofof}{Proposition~\ref{2pnull}}
The sub-scheme $X_{\oes} \subset X_\FS$ is open since $\oes$
is contains $\nu$ whenever $\nu$ is a face of a cone in $\oes$.
Furthermore, $X_{\oes}$ is of finite type over $B$
since $\oes$ is finite. 
Finally, the condition from Remark~\ref{1aG}
is satisfied because any
element of $\oes$ is a face of only finitely many cones.  
\end{Proofof}

Thus, we can consider  
$\Mcogm (X_\FT,\iTjZ) \in \DM$, the 
completed motive of
$X_\FT$ with coefficients in $\iTjZ$ defined in \ref{1aK}. 
The aim of this section is to identify this object. \\

Note that the open covering $(X_{\oes \cap \FT})_{\tau \in \FT}$
of $X_\FT$ corresponds to a closed
covering $(\tau \cap D)_{\tau \in \FT}$ of $D$ with the same
combinatorics. Denote by 
\[
C_\ast \left( (\tau \cap D)_{\tau \in \FT}, \BZ \right)
\]
the homological \v{C}ech complex associated to this covering.
Its $p$-th term
is the direct sum of the $H_0 (\tau_0 \cap \ldots \cap \tau_p \cap D, \BZ)$, 
for all $(p+1)$-fold intersections $\tau_0 \cap \ldots \cap \tau_p$
of the $\tau_i$, 
with $\tau_i \in \FT$. Note that $\tau_0 \cap \ldots \cap \tau_p \cap D$ 
is either contractible of empty according to
whether $\tau_0 \cap \ldots \cap \tau_p$ is a cone in $\FT$ or not.
In particular, the above \v{C}ech complex computes singular homology of 
$D$ with coefficients in $\BZ$.
We define $C_\ast (\FT, \BZ) = (C_p , \delta_p)$, the 
\emph{co-cellular complex of $D$ with respect to the 
stratification induced by $\FT$}, as follows: we set
\[
C_p := \bigoplus_{\tau \in \FT_{-p}} \BZ (\tau) \; ,
\] 
where $\FT_q := \{ \tau \in \FT \tei \dim \tau = q \}$,
and $\BZ (\tau)$ is the \emph{group of orientations} of $\tau$,
i.e., the free Abelian group of rank one defined as the maximal exterior
power
\[
det (\BZ[\nu \in \FS \tei \dim \nu = 1 \; , \; 
\nu \; \text{face of} \; \tau]) \; .
\]
Any ordering of the set $\{ \nu \in \FS \tei \dim \nu = 1 \; , \; 
\nu \; \text{face of} \; \tau \}$ gives a generator of $\BZ (\tau)$.
A change of the ordering by a permutation multiplies the
generator by the signature of the permutation. The tensor square
$\BZ (\tau)^{\otimes 2}$ is canonically isomorphic to $\BZ$.
The restriction of the differential $\delta_p: C_{p+1} \to C_p$ to
$\BZ (\tau)$, for $\tau \in \FT_{- (p+1)}$ is defined as follows:
fix an ordering $\nu_1 \prec \ldots \prec \nu_{- (p+1)}$ of the faces of
dimension one of $\tau$. For $\eta \in \FT_{- p}$, the component
of $\delta_p$ involving $\eta$ is zero if $\tau$ is not a face of $\eta$.
If $\tau$ is a face of $\eta$, then $\eta$ has exactly one more face
$\nu_{-p}$ of dimension one. Define
\[
\BZ (\tau) \longto \BZ (\eta)
\]
by $z \otimes (\nu_1 \prec \ldots \prec \nu_{- (p+1)})
\mapsto z \otimes (\nu_1 \prec \ldots \prec \nu_{- (p+1)} \prec \nu_{-p})$.
This definition is obviously independent of the choice of the ordering
of the set $\{ \nu \in \FS \tei \dim \nu = 1 \; , \; 
\nu \; \text{face of} \; \tau \}$. 
Observe that for any cone
$\eta$ in $\FT$ of dimension $r$ having a face $\tau$ of dimension $r-2$
which still belongs to $\FT$, there are exactly two faces of $\eta$
of dimension $r-1$ having $\tau$ as face, and belonging to $\FT$. This
shows that $C_\ast (\FT, \BZ)$ is indeed a complex.

\begin{Rem} \label{2pBa}
One sees that $C_\ast (\FT', \BZ)$ can be defined for any locally closed
subset $\FT'$ of $\FT$. It is contravariantly functorial with respect
to inclusions $\iota: \FT'' \into \FT'$ of two such subsets:
for $\tau \in \FT'$, the map $\iota^*$ is the identity on $\BZ (\tau)$
if $\tau \in \FT''$, and is zero if $\tau \not \in \FT''$. 
\end{Rem}

Define $D_\BR$ to be the sub-vector space of $Y_*(T)_\BR$ generated by $D$,
and let $\BZ (D_\BR)$ denote the group of orientations of $D_\BR$.
Set $d := \dim_\BR D_\BR$. 
   
\begin{Prop} \label{2pB}
Assume that $D$ is open in $D_\BR$. \\[0.1cm]
(a)~The co-cellular complex of $D$ computes singular homo\-lo\-gy of $D$. 
More precisely, there is a canonical isomorphism
\[
H_p ( C_\ast (\FT, \BZ)) \isoto H_{p+d} (D,\BZ) \otimes_{\BZ \,} \BZ (D_\BR) 
\]
for all $p \in \BZ$. \\[0.1cm]
(b)~If $D$ is contractible, then 
$C_\ast (\FT, \BZ)$ is
a (bounded) resolution of $\BZ (D_\BR)[-d]$.
\end{Prop}

\begin{Proof}
We leave it to the reader to canonically
represent $C_\ast (\FT, \BZ)$ as a sub-quotient of the shift by $-d$ of
$C_\ast \left( (\tau \cap D)_{\tau \in \FT}, \BZ \right) \otimes_{\BZ \,} 
\BZ (D_\BR)$;
note that any choice of orientation of $D_\BR$ induces orientations on
all maximal cones of $\FT$.
The hypothesis on $D$ implies that any cone in $\FT$ of dimension $d-p$, for
$p \ge 0$, is the intersection of precisely $p+1$ maximal cones
in $\FT$. This means that
components $H_0 (\tau_0 \cap \ldots \cap \tau_p \cap D, \BZ)$ of degree $p$ in 
$C_\ast \left( (\tau \cap D)_{\tau \in \FT}, \BZ \right)$, but with
\[
\dim (\tau_0 \cap \ldots \cap \tau_p) < d-p
\]
occur also in degree $p+1$. Using this, one shows that $C_\ast (\FT, \BZ)$
is quasi-isomorphic to the \v{C}ech complex, tensored with
$\BZ (D_\BR)$ and shifted by $- d$.
\end{Proof}
\forget{
One defines $C^\ast (\FT, \BZ)$, the 
\emph{cellular complex of $D$ with respect to the 
stratification induced by $\FT$}, as the dual of $C_\ast (\FT, \BZ)$.
Thus,
\[
C^p  (\FT, \BZ):= \prod_{\tau \in \FT_{- p}} \BZ (\tau) \; .
\] 
This resembles, but should not be confused with the cellular
complex associated to a $CW$-complex. \\
}

Recall the definition of the functor
\[
\Hom: C^b(\Ab) \times \ShN \longto C^b(\ShN)
\]
from Section~\ref{1b}.
Here is the main result of this section:

\begin{Thm} \label{2pC}
(a)~There is a canonical isomorphism
\[
\Mcogm (X_\FT,\iTjZ) \isoto 
\RC \left( \Hom (C_\ast (\FT, \BZ),L(X)) \right) 
\]
in $\DM$. \\[0.1cm]
(b)~If $\FT$ is finite, then
\[
\Mgm (X_\FT,\iTjZ) \isoto 
\RC \left( \Hom (C_\ast (\FT, \BZ),L(X)) \right) 
\]
in $\DM$. \\[0.1cm]
(c)~If $D$ is contractible and open in $D_\BR$, then
\[
\Mcogm (X_\FT,\iTjZ) \isoto \Hom (\BZ (D_\BR),\Mgm (X)) [d]  
\]
in $\DM$. In particular, any choice of an orientation of $D_\BR$
induces an isomorphism
$\Mcogm (X_\FT,\iTjZ) \cong \Mgm (X) [d]$. 
\end{Thm}

Of course, parts~(b) and (c) follow from part~(a), together with
Propositions~\ref{1aL} and \ref{2pB}~(b).
Actually, the proof will show that the above isomorphisms satisfy
strong functoriality properties in the geometric data we fixed in
the beginning. 
For our application to the Baily--Borel compactification of pure
Shimura varieties (Section~\ref{4}), we
need to spell out functoriality in a particular situation: in 
addition to the set of data (1)--(3) and conditions $(\alpha)$--$(\gamma)$
fixed so far, we assume given
(4)~an abstract group $H$ 
satisfying condition (B) from Section~\ref{1b}, 
acting on $B$, $X$, and $T$ in a way compatible
with the group and torsor structures, stabili\-zing $\FS$ and $\FT$,
and satisfying the following condition:
$(\delta)$~the
action of $H$ on the set $\FT$ is free.
This implies that the induced action of $H$ on $X_\FT$ is free and proper. 
As was explained before Theorem~\ref{1bK}, this allows to define
$\Mcogm (X_\FT,\iTjZ) \in H \text{-}\DM$. Observe that the right hand
sides of the isomorphisms from Theorem~\ref{2pC} also define objects
in $H \text{-}\DM$; note that $H$ acts on $\BZ (D_\BR)$
via a character of order at most two. Our proof of Theorem~\ref{2pC} will show:

\begin{Comp} \label{2pCa}
(a)~There is a canonical isomorphism
\[
\Mcogm (X_\FT,\iTjZ) \isoto 
\RC \left( \Hom (C_\ast (\FT, \BZ),L(X)) \right) 
\]
in $H \text{-}\DM$, which is mapped to the one from 
Theorem~\ref{2pC}~(a) under the forgetful functor. \\[0.1cm]
(b)~If $D$ is contractible and open in $D_\BR$, then
\[
\Mcogm (X_\FT,\iTjZ) \isoto \Hom (\BZ (D_\BR),\Mgm (X)) [d]  
\]
in $H \text{-}\DM$. This isomorphism is compatible with the isomorphism from
Theorem~\ref{2pC}~(c).
\end{Comp}

Note that under the above conditions,
the action of $H$ is also free and proper
on the formal completions $(X_{\FS - \{ 0 \}})_{X_\FT}$ 
and $(X_\FS)_{X_\FT}$ of $X_{\FS - \{ 0 \}}$ and $X_\FS$ along
$X_\FT$. Hence the geometric quotients $H \backslash X_\FT$,
$H \backslash (X_{\FS - \{ 0 \}})_{X_\FT}$, and $H \backslash (X_\FS)_{X_\FT}$
exist, and $X_\FT \to H \backslash X_\FT$,
$(X_{\FS - \{ 0 \}})_{X_\FT} \to H \backslash (X_{\FS - \{ 0 \}})_{X_\FT}$, and
$(X_\FS)_{X_\FT} \to H \backslash (X_\FS)_{X_\FT}$
are local isomorphisms. \\

We make one last extension of the data:
we assume given (5)~separated reduced schemes $S$,
$S_\FS$, and $S_\FT$, which are locally of finite type over $k$, and
related by an open immersion $S \into S_\FS$ and a closed immersion
$S_\FT \into S_\FS$ identifying $S_\FT$ with a closed sub-scheme of
the reduced scheme $S_{\FS - \{ 0 \}}:= S_\FS - S$. 
We assume that $S_\FT \in \schi$. 
By abuse of notation, we again use the symbols
$j$ and $i_\FT$ for these immersions. We fix 
(6)~a triple of compatible isomorphisms
\[ 
H \backslash (X_\FS)_{X_\FT} \isoto (S_\FS)_{S_\FT}
\; , \; 
H \backslash (X_{\FS - \{ 0 \}})_{X_\FT} \isoto (S_{\FS - \{ 0 \}})_{S_\FT}
\; , \; 
H \backslash X_\FT \isoto S_\FT \; ,
\]
all denoted by the same symbol $f$,
where $(\argdot)_{S_\FT}$ denotes the formal completions along $S_\FT$.
Complement~\ref{2pCa} and Theorem~\ref{1bK} then imply:

\begin{Thm} \label{2main}
(a)~Under the above hypotheses, $f$ induces an isomorphism
\[
\Mcogm (S_\FT,\iTjZ) \isoto 
\RC \left( \Hom (C_\ast (\FT, \BZ),L(X))^H \right)
\]
in $\DM$. \\[0.1cm]
(b)~If $H \backslash \FT$ is finite (hence $S_\FT \in Sch/k$), then
\[
\Mgm (S_\FT,\iTjZ) \isoto 
\RC \left( \Hom (C_\ast (\FT, \BZ),L(X))^H \right) 
\]
in $\DM$. \\[0.1cm]
(c)~If $D$ is contractible and open in $D_\BR$, then
\[
\Mcogm (S_\FT,\iTjZ) \isoto 
R \Gamma \left(H, \Hom (\BZ (D_\BR),\Mgm (X)) \right)[d]  
\]
in $\DM$. 
\end{Thm}

\begin{Proof}
The isomorphisms in $\DM$ are those of Theorem~\ref{1bK}.
That the complexes in (a) resp.\ (b) represent the right classes
in $\DM$ results from Corollary~\ref{1bD}~(b) and hypothesis $(\delta)$
on the action of $H$ on $\FT$. 
\end{Proof}

Note that hypotheses (b) and (c) of Theorem~\ref{2main} together imply
that $H$ is of type $FL$. \\

We prepare the proof of Theorem~\ref{2pC}
(and forget the data (4)--(6) and condition $(\delta)$). We need to
identify the contribution of each $X_{\os}$
to the completed motive $\Mcogm (X_\FT,\iTjZ)$.  
Write $n := \dim T$ and fix a cone $\sigma \in \FS$
of dimension $r$. 
The finite subset $\os$ of $\FS$ corresponds to an open sub-scheme
$X_{\os} \subset X_\FS$ which inherits the stratification,
and which is of finite type over $B$.
For any face $\nu$ of $\sigma$, the scheme $X_{\ot}$
is open in $X_{\os}$, and contains $X_\nu$ as unique closed stratum.
Use the symbol $j$ to denote also the open immersion of $X$ into $X_{\os}$.
Zariski-locally over $B$, $j: X \into X_{\os}$ is
isomorphic to $\Gnm \into \Gnrrm \times \BA^r$, and
the stra\-ti\-fication is the natural strati\-fi\-cation of $\BA^n$
according to the vanishing, resp.\ non-vanishing of coordinates. 
Let us start with the identification of
$\Mgm (\Xs,\isjZ)$. By definition,
this is the motive $\Mgm (\FY_\sigma)$ associated to the diagram 
\[
\xymatrix@R-20pt{
& X_{\os - \{ 0, \sigma \}} \ar[r] \ar[dd] & 
X_{\os - \{ \sigma \}} \ar[dd] \\
\FY_\sigma \quad := & & \\
& X_{\os - \{ 0 \}}  \ar[r] & 
X_{\os} 
\\}
\] 
as in \cite[Conv.~1.2]{W9}. The component $X_{\os}$ sits in total degree
zero. Our main computational 
tool is the following:

\begin{Lemma} \label{2A}
Let $\nu \in \os$. The immersion
$X_{\ot - \{ 0 \}} \into X_{\ot}$ induces an isomorphism
\[
\Mgm (X_{\ot - \{ 0 \}}) \isoto \Mgm(X_{\ot}) 
\]
in $\DM$ provided that $\nu \ne 0$. 
\end{Lemma}

\begin{Proof}
The Mayer--Vietoris property
\cite[Prop.~3.1.3]{VSF} for the functor $L$
shows that we may assume that
$X = T \times_k B$. For a suitable choice of
identification $\Gnm \cong T$, the immersion $X_0 \into X_{\ot}$
is isomorphic to 
\[
\Gnm \longinto \Gnrm \times \BA^i \; , 
\]
with $i \ne 0$. We may reduce to the case $i = n$. Now observe that
\[
\BA^n - \Gnm \longinto \BA^n 
\]
is a homotopy equivalence for $n \ne 0$ (in fact, both sides are
homotopic to a point).
\end{Proof}

In order to compute the motive $\Mgm (\FY_\sigma)$, 
consider the upper line 
\[
X_{\os - \{ 0, \sigma \}} \longto
X_{\os - \{ \sigma \}} 
\]
of $\FY_\sigma$, and the finite covering 
\[
X_{\os - \{ \sigma \}} = \bigcup_{\nu \in \os^{\, 1}} X_{\ot} \; ,
\]
where $\os^{\, 1} \subset \os$ denotes the subset (of cardinality $r$)
of faces $\nu$
of codimension one in $\sigma$. For $\nu_1, \ldots, \nu_p \in \os^{\, 1}$,
we have:
\[
\bigcap_{i=1}^p {X_{\oit}} = X_{\ot} \; , \; 
\text{ with } \; \nu = \cap_i \nu_i \; .
\]
Fix an ordering of $\{ \nu \in \FS \tei \dim \nu = 1 \; , \; 
\nu \; \text{face of} \; \sigma \}$. This induces an ordering of
$\os^{\, 1}$.
Define the double complex $L(\FY'_\sigma)$ of Nisnevich sheaves 
as follows: first consider the \v{C}ech complex associated to the
above covering, and to the induced covering of 
$X_{\os - \{ 0, \sigma \}}$, then add the lower line of $L(\FY_\sigma)$: 
\[
\xymatrix@R-20pt{
0 = L(\emptyset) \ar[r] \ar[dd] &
L(X) = L(X_0) \ar[dd] \\
& \\
\vdots \ar[dd] &
\vdots \ar[dd] \\
& \\
\prod_{\nu_1,\nu_2 \in \os^{\, 1} \atop \nu_1 \prec \nu_2} 
L(X_{\overline{ \{ \nu_1 \cap \nu_2 \}} - \{ 0 \}}) \ar[r] \ar[dd] &
\prod_{\nu_1,\nu_2 \in \os^{\, 1} \atop \nu_1 \prec \nu_2} 
L(X_{\overline{ \{ \nu_1 \cap \nu_2 \} }}) \ar[dd] \\
& \\
\prod_{\nu \in \os^{\, 1}} L(X_{\ot - \{ 0 \}}) \ar[r] \ar[dd] & 
\prod_{\nu \in \os^{\, 1}} L(X_{\ot}) \ar[dd] \\
& \\
L(X_{\os - \{ 0 \}})  \ar[r] & 
L(X_{\os}) 
\\}
\] 
Again, we have total
degree zero for $L(X_{\os})$, hence degree $-r$ for $L(X)$.
We get:

\begin{Lemma} \label{2B}
Fix an ordering of 
$\{ \nu \in \FS \tei \dim \nu = 1 \; , \; 
\nu \; \text{face of} \; \sigma \}$. \\[0.1cm]
(a)~There is a canonical morphism of complexes of Nisnevich sheaves
\[
MV_\sigma: L(\FY'_\sigma) \longto L(\FY_\sigma) \; .
\]
It is a quasi-isomorphism, and hence
induces an isomorphism
$\RC (MV_\sigma)$ in $\DM$. \\[0.1cm]
(b)~The projection from $L(\FY'_\sigma)$ to its upper line 
\[
pr_\sigma: L(\FY'_\sigma) \longto L(X) [\dim \sigma]
\]
induces an isomorphism
$\RC (pr_\sigma)$ in $\DM$. 
\end{Lemma} 

\begin{Proof}
(a) is the Mayer--Vietoris property, and (b)
is a consequence of Lemma~\ref{2A}.
\end{Proof}

Analyzing the dependence of the isomorphisms on the choice of the
ordering, one finds:

\begin{Cor} \label{2C}
Let $\sigma \in \FS$. There is a canonical isomorphism
\[
\Mgm (\Xs,\isjZ) \isoto \Hom (\BZ (\sigma),\Mgm (X))[\dim \sigma] 
\]
in $\DM$. Any choice of an ordering $\prec$ as in Lemma~\ref{2B}
induces an isomorphism $\rho_\prec$ of the
right hand side with $\Mgm (X)[\dim \sigma]$. The induced isomorphism
\[
\Mgm (\Xs,\isjZ) \isoto \Mgm (X)[\dim \sigma] 
\]
equals $\RC (pr_\sigma) \circ (\RC (MV_\sigma))^{-1}$, where 
$pr_\sigma$ and $MV_\sigma$ are as in Lemma~\ref{2B}.
\end{Cor}

We are ready for the proof of Theorem~\ref{2pC}~(a)
in the affine case, i.e.,
with $\FS$ replaced by $\oes$, for $\tau \in \FT$: 
recall that $\Mgm (X_{\FT \cap \oes},\iTjZ)$
is associated to  
\[
\xymatrix@R-20pt{
& X_{\oes - \FT - \{ 0 \}} \ar[r] \ar[dd] & 
X_{\oes - \FT} \ar[dd] \\
\FY_{\FT \cap \oes} \quad := & & \\
& X_{\oes - \{ 0 \}}  \ar[r] & 
X_{\oes} 
\\}
\] 

\begin{Lemma} \label{2nA}
Fix an ordering $\prec$ of 
$\{ \nu \in \FS \tei \dim \nu = 1 \; , \; 
\nu \; \text{face of} \; \tau \}$. \\[0.1cm]
(a)~There is a complex of Nisnevich sheaves $L(\FY'_{\FT \cap \oes})$, 
and a canonical quasi-isomorphism
\[
MV_{\FT \cap \oes}: L(\FY'_{\FT \cap \oes}) \longto L(\FY_{\FT \cap \oes}) \; .
\]
(b)~There is a morphism 
\[
\rho_\prec^{-1} \circ pr_{\FT \cap \oes}: L(\FY'_{\FT \cap \oes}) \longto 
\Hom (C_\ast (\FT \cap \oes, \BZ),L(X))
\]
inducing an isomorphism in $\DM$. \\[0.1cm] 
(c)~The resulting isomorphism
\[
\Mgm (X_{\FT \cap \oes},\iTjZ) \isoto 
\RC \left( \Hom (C_\ast (\FT \cap \oes, \BZ),L(X)) \right) 
\]
in $\DM$ does not depend on the choice of the ordering. \\[0.1cm]
(d)~If $\tau$ is a face of $\sigma \in \FT$, then any extension of $\prec$
to an ordering of $\{ \nu \in \FS \tei \dim \nu = 1 \; , \; 
\nu \; \text{face of} \; \sigma \}$ induces a canonical morphism of
complexes
\[
L(\FY'_{\FT \cap \oes}) \longto L(\FY'_{\FT \cap \os})
\] 
such that
\[
\xymatrix@R-20pt{
L(\FY'_{\FT \cap \oes}) \ar[rr]^-{MV_{\FT \cap \oes}} \ar[dd] & & 
L(\FY_{\FT \cap \oes}) \ar[dd] \\
& & \\
L(\FY'_{\FT \cap \os}) \ar[rr]^-{MV_{\FT \cap \os}} & &  
L(\FY_{\FT \cap \os})
\\}
\] 
and
\[
\xymatrix@R-20pt{
L(\FY'_{\FT \cap \oes}) \ar[rr]^-{\rho_\prec^{-1} \circ pr_{\FT \cap \oes}} 
                                                                   \ar[dd] & & 
\Hom (C_\ast (\FT \cap \oes, \BZ),L(X)) \ar[dd] \\
& & \\
L(\FY'_{\FT \cap \os}) \ar[rr]^-{\rho_\prec^{-1} \circ pr_{\FT \cap \os}} & &  
\Hom (C_\ast (\FT \cap \os, \BZ),L(X)) \\
\\}
\] 
commute. Here, the morphism $\FY_{\FT \cap \oes} \longto \FY_{\FT \cap \os}$
is the natural one given by the inclusion $\iota$ of $\oes$ into 
$\os$, and 
\[
\Hom (C_\ast (\FT \cap \oes, \BZ),L(X)) \longto
\Hom (C_\ast (\FT \cap \os, \BZ),L(X))
\]
is induced by the morphism $\iota^*$ from Remark~\ref{2pBa}.
\end{Lemma}

\begin{Proof}
The lower row of $\FY_{\FT \cap \oes}$
contributes nothing to the motive $\Mgm (X_{\FT \cap \oes},\iTjZ)$
thanks to Lemma~\ref{2A} (note that $\tau \ne 0$
since $\tau \in \FT$). For the
upper row
\[
\FZ_{\FT \cap \oes} \quad := \quad
X_{\oes - \FT - \{ 0 \}} \longto X_{\oes - \FT} \; ,
\]
we shall construct a double complex $L(\FZ'_{\FT \cap \oes})$
of vertical length $\dim \tau$, and
resolving $L(\FZ_{\FT \cap \oes})$ as in the construction 
preceding Lemma~\ref{2B}. 
Every row in $L(\FZ'_{\FT \cap \oes})$ will be the direct product
of one of two types of factors: either, the factor is simply
\[
0 = L(\emptyset) \longto L(X) = L(X_0) \; ,
\]
or it has the shape
\[
L(X_{\oae - \{ 0 \}})  \longto L(X_{\oae}) \; ,
\]
for a face $\eta$ of $\tau$ not belonging to $\FT$, and unequal to $0$. 
The vertical differentials of $L(\FZ'_{\FT \cap \oes})$ will be
induced by the inclusions of faces.

For the construction of $L(\FZ'_{\FT \cap \oes})$,
we apply induction on the cardinality of $\FT \cap \oes$,
starting with Lemma~\ref{2B} and Corollary~\ref{2C} if
$\FT \cap \oes = \{ \tau \}$. In the general case,
consider the finite covering
\[
X_{\oes - \FT} = \bigcup_{\nu \in \oes^{\, 1}} X_{\ot - \FT} \; ,
\]
and observe that for each $\nu \in \oes^{\, 1}$, the cone $\tau$ is not
contained in the set
$\FT \cap \ot$. Hence the latter is properly contained in $\FT \cap \oes$,
and we may apply the induction step on each 
\[
X_{\oae - \FT - \{ 0 \}} \longto X_{\oae - \FT} \; ,
\]
for any intersection $\eta$ of cones in $\oes^{\, 1}$
contained in $\FT$. Note that for $\eta \not \in \FT$, this
row equals
\[
X_{\oae - \{ 0 \}} \longto X_{\oae} \; .
\]
The result is the double complex $L(\FZ'_{\FT \cap \oes})$.
Note that the top row of $L(\FZ'_{\FT \cap \oes})$ equals (one copy of)
\[
0 \longto L(X) \; .
\]
By definition, we get
$L(\FY'_{\FT \cap \oes})$ by adding the lower row of $L(\FY_{\FT \cap \oes})$.
The canonical map from $L(\FY'_{\FT \cap \oes})$ 
to $L(\FY_{\FT \cap \oes})$
is a quasi-isomor\-phism thanks to the Mayer--Vietoris property.
The factors of the shape 
\[
L(X_{\oae - \{ 0 \}}) \longto L(X_{\oae}) \; ,
\]
with $\eta \ne 0$, organize into a sub-complex of $L(\FY'_{\FT \cap \oes})$.
(Note: each such $\eta$ is either a face of $\tau$ not belonging to
$\FT$, or equal to $\tau$, in which case the factor in question
is the last row of $L(\FY'_{\FT \cap \oes})$.)
By Lemma~\ref{2A}, this sub-complex vanishes in $\DM$. 
Furthermore, by induction, the quotient of $L(\FY'_{\FT \cap \oes})$
by this sub-complex consists of the top row 
\[
0 \longto L(X) 
\]
of $L(\FY'_{\FT \cap \oes})$, and the 
$\Hom (C_\ast (\FT \cap \oae, \BZ),L(X))$,
for all intersections $\eta$ of cones in $\oes^{\, 1}$
contained in $\FT$.
More precisely, this quotient equals
\[
\Hom (C'_\ast (\FT \cap \oes, \BZ),L(X)) \; ,
\]
where $C'_\ast (\FT \cap \oes, \BZ)$ is the complex defined as follows:
for $m > -\dim \tau$, we have
\[
C'_m (\FT \cap \oes, \BZ) :=
\bigoplus_{p+q = m} \bigoplus_\eta C_p (\FT \cap \oae, \BZ) \; ,
\]
where $\eta$ runs through all $(-q+1)$-fold intersections
of cones in $\oes^{\, 1}$. We have
\[
C'_{-\dim \tau} (\FT \cap \oes, \BZ) := \BZ(\tau) \; ,
\]
and $C'_m (\FT \cap \oes, \BZ) := 0$ for 
$m \not \in [-\dim \tau , 0]$. The differential
\[
C'_{-(\dim \tau - 1)} (\FT \cap \oes, \BZ)
\longto C'_{-\dim \tau} (\FT \cap \oes, \BZ)
\]
is the direct sum of maps $\BZ(\eta) \to \BZ(\tau)$,
for faces $\eta$ of $\tau$. These maps are all isomorphisms.
On orientations, they are described 
by adding the missing faces of dimension one,
in the fixed order. It remains to construct a morphism
\[
i: C_\ast (\FT \cap \oes, \BZ) \longto C'_\ast (\FT \cap \oes, \BZ) \: .
\]
In degree $-\dim \tau$, we take the identity on $\BZ(\tau)$.
In degree $m > -\dim \tau$, 
take $\sigma \in \FT \cap \oes$ of
dimension $-m$, so $\BZ (\sigma)$ is a summand of $C_m (\FT \cap \oes, \BZ)$.
It is also a summand of 
$C_m (\FT \cap \oae, \BZ)$, 
for any face $\eta$ of $\tau$ containing $\sigma$ as a face,
and we define $i_m$ on $\BZ (\sigma)$
as the diagonal embedding of $\BZ (\sigma)$ in 
$\oplus_\eta C_m (\FT \cap \oae, \BZ) \subset C'_m (\FT \cap \oes, \BZ)$.
We leave it to the reader to check that this defines 
a morphism of complexes, and in fact, a quasi-isomorphism.
Since both $C_* (\FT \cap \oes, \BZ)$
and $C'_* (\FT \cap \oes, \BZ)$ are bounded
complexes of free Abelian groups,
it follows that $i$ is even a homotopy equivalence. Hence the
induced morphism
\[
\Hom (C'_* (\FT \cap \oae, \BZ),L(X)) 
\longto \Hom (C_\ast (\FT \cap \oae, \BZ),L(X)) 
\]
is a homotopy equivalence as well.
\end{Proof}

\forget{
It should be remarked that in the geometrical setting of Lemma~\ref{2nA}~(d),
but without the compatibility condition on the orderings, it is still true
that there is a canonical morphism of complexes
\[
L(\FY'_{\FT \cap \oes}) \longto L(\FY'_{\FT \cap \os})
\] 
making the two diagrams commutative. \\
}
\begin{Proofof}{Theorem~\ref{2pC}}
By Proposition~\ref{2pnull}~(b), the
system $\FU_{X_\FS}$ given as $(X_{\oes})_{\tau \in \FT}$ 
lies in $\Cov (X_\FS)_{X_\FT}$.
By Definition~\ref{1aK} and the remark following it, the completed motive
$\Mcogm (X_\FT,\iTjZ)$ is represented by the simple complex
associated to $\uC_* (L^\bullet (\FU_{\FY_\FT}))$, with
\[
\xymatrix@R-20pt{
& X_{\FS - \FT - \{ 0 \}} \ar[r] \ar[dd] & 
X_{\FS - \FT} \ar[dd] \\
\FY_\FT \quad := & & \\
& X_{\FS - \{ 0 \}}  \ar[r] & 
X_\FS 
\\}
\] 
Recall that $L^\bullet (\FU_{\FY_\FT})$ is the simple complex
associated to
\[
\xymatrix@R-20pt{
L^\bullet (\FU_{X_{\FS - \FT - \{ 0 \}}}) \ar[r] \ar[dd] &
L^\bullet (\FU_{X_{\FS - \FT}}) \ar[dd] \\
& \\
L^\bullet (\FU_{X_{\FS - \{ 0 \}}}) \ar[r] &
L^\bullet (\FU_{X_{\FS}})
\\}
\]
with $L^\bullet (\FU_{X_{\FS}})$ sitting in degree $(0,\argast)$.
We may represent $L^\bullet (\FU_{\FY_\FT})$ as follows: define
a triple complex whose term $(-m,\argast,\argast)$
is the direct product of the $L(\FY_{\FT \cap \oes})$, 
for all $(m+1)$-fold intersections $\tau$ in $\FT$.
Then $L^\bullet (\FU_{\FY_\FT})$ is the simple complex associated
to this triple complex.
By Lemma~\ref{2nA}, $L^\bullet (\FU_{\FY_\FT})$ may be replaced
by the simple complex $s L^{\bullet,\bullet}$ 
associated to $L^{\bullet,\bullet}$, with
\[
L^{p,q} := \prod_\tau \Hom (C_p (\FT \cap \oes, \BZ),L(X)) \; ,
\]
where $\tau$ runs through all $(-q+1)$-fold intersections in $\FT$.
Note that $L^{p,q}$ is non-zero only if $p$ and $q$ are non-positive.

We need to define a morphism
\[
s L^{\bullet,\bullet} \longto \Hom (C_\ast (\FT, \BZ),L(X)) \; ,
\]
or in other words, a family of morphisms
\[
\prod_{p+q = m} L^{p,q} \longto \Hom (C_m (\FT, \BZ),L(X))
\]
compatible with the differentials. We shall define the dual morphisms
\[
i_m: C_m (\FT, \BZ) \longto 
\bigoplus_{p+q = m} \bigoplus_\tau C_p (\FT \cap \oes, \BZ) \; :
\]
Let $\eta \in \FT$ be of
dimension $-m$, so $\BZ (\eta)$ is a summand of $C_m (\FT, \BZ)$.
It is also a summand of 
$C_m (\FT \cap \oes, \BZ)$, 
for any $\tau$ containing $\eta$ as a face,
and we define $i_m$ on $\BZ (\eta)$
as the diagonal embedding of $\BZ (\eta)$ in 
$\oplus_\tau C_m (\FT \cap \oes, \BZ)$, where the direct sum extends over
all (finitely many!) $\tau$ containing $\eta$ as a face.
This defines the morphism of complexes
\[
i: C_* (\FT, \BZ) \longto s C_{*,*} (\FT, \BZ) \; ,
\]
where the right hand side denotes the simple complex associated to
$C_{*,*} (\FT, \BZ)$, with 
\[
C_{p,q} (\FT, \BZ) := \bigoplus_\tau C_p (\FT \cap \oes, \BZ) 
\]
as above. One sees that $i$ factors through
$Z_{*,0} (\FT, \BZ) \subset s C_{*,*} (\FT, \BZ)$, with
$Z_{p,q} (\FT, \BZ)$ 
denoting the kernel of the differential 
in $q$-direction
\[
d_q: C_{p,q} (\FT, \BZ) \longto C_{p,q-1} (\FT, \BZ) \: .
\]
Given the definition of that differential, it is not difficult to check
that the complex
\[
C_m (\FT, \BZ) \stackrel{i_m}{\longto} C_{m,0} (\FT, \BZ) 
\stackrel{d_0}{\longto} C_{m,-1} (\FT, \BZ) 
\stackrel{d_{-1}}{\longto} \ldots
\]
is acyclic for every $m$. The complex $C_* (\FT, \BZ)$ is bounded,
hence $i$ is a quasi-isomorphism. Since both $C_* (\FT, \BZ)$
and $sC_{*,*} (\FT, \BZ)$ are complexes of free Abelian groups,
and $C_* (\FT, \BZ)$ is bounded,
it follows that $i$ is even a homotopy equivalence. Hence the
induced morphism
\[
s L^{\bullet,\bullet} 
= \Hom (sC_{*,*} (\FT, \BZ),L(X)) 
\longto \Hom (C_\ast (\FT, \BZ),L(X)) 
\]
is a homotopy equivalence as well.
\forget{
Consider the ascending filtration on $\FT$ defined by
\[
\FT_{\ge -m} := \{ \tau \in \FT \tei \dim \tau \ge -m \} \; .
\]
It induces filtrations by sub-complexes
on both the target and the source of 
$i$. Furthermore, by definition, the morphism $i$ respects 
these filtrations.
In order to prove our claim, it suffices to show that the 
morphisms induced on the graded parts  
are quasi-isomorphisms.

As for the source of $i$, we get the complex concentrated
in degree $m$ and with value $C_m (\FT, \BZ)$.

Variant~\ref{1aN} implies that the graded parts of the source compute the 
completed motives
$\Mcogm (X_{\FT_{n-m}},\iTnmjZ)$,
where $i_{\FT_{n-m}}$ is the immersion of the sub-scheme
$X_{\FT_{n-m}}$. The graded parts of the target are of the simple shape
\[
\prod_{\dim \tau = n-m} \Hom (\BZ (\tau),L(X))[n-m] \; .
\]
For each $m$, the morphism
\[
\Mcogm (X_{\FT_{n-m}},\iTnmjZ) \longto 
\prod_{\tau \in \FT_{n-m}} \Hom (\BZ (\tau),L(X))[n-m] 
\]
is the one coming from the system $(X_{\oes})_{\tau \in \FT_{n-m}}$ 
in $\Cov (X_\FS)_{X_{\FT_{n-m}}}$. By Corollary~\ref{2C}, it is 
an isomorphism in $\DM$.
}
\end{Proofof}

\begin{Ex}
Note that the conclusion of Theorem~\ref{2pC}~(c) is false in general
if $D$ is only contractible, but not open. As an example, take
$\FT$ to consist of two cones: $\sigma$, of dimension at least two,
and (only) one of the faces of $\sigma$ of codimension one. 
The conclusion of \ref{2pC}~(c) is 
false since $\Mgm (X_\FT,\iTjZ) = 0$ by \ref{2pC}~(b).
\end{Ex}


\bigskip

%
%

\section{Notations for Shimura varieties}
\label{3pre}



We recall the notation for Shimura varieties and their compactifications.
It is identical to the one used in \cite{P1,P2,W1,BW}, except that we change
the capital letter $M$ used to denote Shimura varieties
in \loccit \ to capital $S$ in order to avoid
confusion with the motivic notation introduced earlier in this paper. \\

Let $(P, \FX)$ be \emph{mixed Shimura data} \cite[Def.~2.1]{P1}. 
In particular, $P$ is a connected algebraic linear group over $\BQ$, and
$P(\BR)$ acts on the complex ma\-ni\-fold $\FX$ by analytic automorphisms.
Denote by $W$ the unipotent radical of $P$. If $P$ is reductive, i.e., 
if $W=0$, then $(P, \FX)$ is called \emph{pure}.
The \emph{Shimura varieties} associated to $(P,\FX)$ are indexed by the 
open compact subgroups of $P (\BA_f)$.
If $K$ is one such, then the analytic space of $\BC$-valued points of the
corresponding variety is given as
\[
S^K (\BC) := P (\BQ) \backslash ( \FX \times P (\BA_f) / K ) \; ,
\]
where $\BA_f$ denotes the ring of finite ad\`eles over $\BQ$.
According to
Pink's gene\-ral\-ization to mixed Shimura varieties of the Algebraization
Theorem of Baily and Borel
\cite[Prop.~9.24]{P1}, there exist canonical structures of normal algebraic
varieties on the $S^K(\BC)$, which we denote as
\[
S^K_{\BC} := S^K (P, \FX)_{\BC} \; .
\]
According to
\cite[Thm.~2.18]{M} and \cite[Thm.~11.18]{P1}, there is a \emph{canonical model}
of $S^K_{\BC}$, which we denote as
\[
S^K := S^K (P, \FX) \; .
\]
It is defined over the \emph{reflex field} $E (P , \FX)$ of $(P , \FX)$ 
\cite[11.1]{P1}. \\

Any \emph{admissible parabolic subgroup} \cite[Def.~4.5]{P1} $Q$ of $P$
has a canonical normal subgroup $P_1$ \cite[4.7]{P1}.
There is a finite collection of \emph{rational boun\-dary components} 
$(P_1 , \FX_1)$, indexed by the $P_1 (\BR)$-orbits in $\pi_0 (\FX)$ 
\cite[4.11]{P1}. The $(P_1 , \FX_1)$ are themselves mixed Shimura data. 
Consider the following condition on $(P, \FX)$:
\begin{enumerate}
\item [$(+)$] If $G$ denotes the maximal reductive quotient of $P$, 
then the neutral connected component $Z (G)^0$ of the center $Z (G)$ of 
$G$ is, up to isogeny, a direct product of a $\BQ$-split torus with a torus 
$T$ of compact type (i.e., $T(\BR)$ is compact) defined over $\BQ$.
\end{enumerate}
If $(P, \FX)$ satisfies $(+)$, 
then so does any rational boundary component $(P_1 , \FX_1)$
\cite[proof of Cor.~4.10]{P1}.
Denote by $U_1 \subset P_1$ the ``weight $-2$'' part of $P_1$.
It is Abelian, normal in $Q$, and central in the unipotent radical 
$W_1$ of $P_1$.
Fix a connected component $\FX^0$ of $\FX$, and denote by
$(P_1 , \FX_1)$ the associated rational boundary component. 
There is a natural open embedding
\[
\iota: \FX^0 \longinto \FX_1
\]
\cite[4.11, Prop.~4.15~(a)]{P1}. 
If $\FX^0_1$ denotes the connected component of $\FX_1$ containing $\FX^0$, 
then the image of the embedding can be described by means of the map 
\emph{imaginary part}
\[
\imm : \FX_1 \longrightarrow U_1 (\BR) (-1) := 
\frac{1}{2 \pi i} \cdot U_1 (\BR) \subset U_1 (\BC)
\]
of \cite[4.14]{P1}: $\iota (\FX^0)$ is the pre-image of an open convex cone
\[
C (\FX^0 , P_1) \subset U_1 (\BR) (-1)
\]
under $\imm |_{\FX^0_1}$ \cite[Prop.~4.15~(b)]{P1}. 
The reflex field does not change when passing from $(P , \FX)$ to a 
rational boundary component \cite[Prop.~12.1]{P1}.


\bigskip

%
%

\section{Toroidal compactifications of Shimura varieties}
\label{3}



In order to discuss toroidal compactifications, we need to introduce 
the \emph{co\-ni\-cal complex} associated to $(P,\FX)$: 
set-theoretically, it is defined as
\[
\CC (P , \FX) := \coprod_{(\FX^0 , P_1)} C (\FX^0 , P_1) \; .
\]
By \cite[4.24]{P1}, the conical complex is naturally equipped with a 
topology (usually different from the coproduct topology). 
The closure $C^{\ast} (\FX^0 , P_1)$ of $C (\FX^0 , P_1)$ 
inside $\CC (P , \FX)$
can still be considered as a convex cone in $U_1 (\BR) (-1)$, 
with the induced topology. \\

For a fixed open compact subgroup $K \subset P (\BA_f)$,
the (partial) \emph{toroidal compactifications} of 
$S^K$ are para\-me\-terized by \emph{$K$-admissible partial cone 
decompositions}, which are collections of subsets of
\[
\CC (P,\FX) \times P (\BA_f)
\]
satisfying the axioms of \cite[6.4]{P1}. 
If $\FS$ is one such, then in particular any member of $\FS$ is of the shape
\[
\sigma \times \{ p \} \; ,
\]
$p \in P (\BA_f)$, $\sigma \subset C^{\ast} (\FX^0 , P_1)$ a
convex rational polyhedral cone in 
the vector space $U_1 (\BR) (-1)$ 
not containing any non-trivial linear subspace. \\

Let $S^K (\FS) := S^K (P , \FX , \FS)$ be the associated compactification;
we refer to \cite[9.27, 9.28]{P1}
for criteria sufficient to guarantee its existence. 
It comes equipped with a natural stratification into locally closed strata.
The unique open stratum is $S^K$. Any stratum different from the generic
one is obtained as follows: 
Fix a pair $(\FX^0 , P_1)$ as above, $p \in P (\BA_f)$ and
\[
\sigma \times \{ p \} \in \FS
\]
such that $\sigma \subset C^{\ast} (\FX^0 , P_1)$,
$\sigma \cap C (\FX^0 , P_1) \neq \emptyset$,
and the admissible parabolic subgroup $Q$ giving rise to $P_1$ is \emph{proper},
i.e., unequal to $P$.
To $\sigma$, one associates Shimura data
\[
\left( \Pes , \Xes \right)
\]
\cite[7.1]{P1}, whose underlying group 
$\Pes$ is the quotient of $P_1$ by the algebraic subgroup
\[
\langle \sigma \rangle \subset U_1
\]
satisfying 
$\BR \cdot \sigma = \frac{1}{2 \pi i} \cdot \langle \sigma \rangle (\BR)$. Set 
\[
K_1  :=  P_1 (\BA_f) \cap p \cdot K \cdot p^{-1} \; , \quad
\pi_{[\sigma]}  :  P_1 \longonto \Pes \; .
\]
According to \cite[7.3]{P1}, there is a canonical map
\[
\isp (\BC) : \Sps (\Pes , \Xes) (\BC) \longrightarrow 
S^K (\FS) (\BC) := S^K (P , \FX , \FS) (\BC)
\]
whose image is locally closed. 
In fact, the source and the target of $\isp (\BC)$ have canonical
models over the reflex field $E (P , \FX)$, and $\isp (\BC)$
comes from a morphism of schemes over $E (P , \FX)$
\cite[Thm.~12.4]{P1} denoted $\isp$.
If $(P, \FX)$ satisfies $(+)$, and $K$ is \emph{neat} 
(see e.g.\ \cite[0.6]{P1}), then $\isp$ is an immersion, i.e., it 
identifies $\Sps$ with a locally closed sub-scheme of $M^K (\FS)$
\cite[Prop.~1.6]{W1}.
Denote by 
\[
j: S^K \longinto S^K (\FS) 
\]
the open immersion, and by $S^{K_1}$ the Shimura variety
$S^{K_1} (P_1 , \FX_1)$.
Recall the definition of $\Mgm (\Sps, i_{\sigma , p}^! \, j_! \, \BZ)$,
the motive of $\Sps$ with coefficients in $i_{\sigma , p}^! \, j_! \, \BZ$
\cite[Def.~3.1]{W9}, and of
the group of orientations 
$\BZ (\sigma)$ of $\sigma$ (Section~\ref{2}).
Here is 
our main result concerning strata in toroidal compactifications:

\begin{Thm} \label{3main}
Assume that $(P, \FX)$ satisfies $(+)$, and that $K$ is neat. 
Then there is a canonical isomorphism 
\[
\Mgm (\Sps, i_{\sigma , p}^! \, j_! \, \BZ) \isoto
\Hom (\BZ (\sigma),\Mgm (S^{K_1})) [\dim \sigma]  
\]
in $DM^{eff}_-(E (P , \FX))$. 
\end{Thm}

\begin{Proof}
Let $\Ses$ be the minimal $K_1$-admissible cone decomposition of
\[
\CC (P_1 , \FX_1) \times P_1 (\BA_f)
\]
containing $\sigma \times \{ 1 \}$. It
can be realized inside the decomposition $\FS^0_1$ of 
\cite[6.13]{P1}; by definition, $\Ses$ is 
\emph{concentrated in the unipotent fibre} \cite[6.5~(d)]{P1}. 
View $\Sps$ as sitting inside $S^{K_1} (\Ses)$:
$i_1: \Sps \into S^{K_1} (\Ses)$.
In fact, $j_1: S^{K_1} \into S^{K_1} (\Ses)$
is a relative affine torus embedding over $\Sps$,
and $i_1$ identifies
$\Sps$ with the unique closed stratum 
of the canonical stratification \cite[Prop.~1.15, Lemma~1.16]{W1}.
Consider the diagram
\[
\vcenter{\xymatrix@R-10pt{ 
\Sps \ar@{^{ (}->}[r]^{\isp} \ar@{_{ (}->}[dr]_{i_1} &
S^K (\FS) \\
&  S^{K_1} (\Ses) \\}}
\]
According to \cite[Thm.~1.13]{W1}, there is a canonical isomorphism 
\[
f: (S^K (\FS))_{\Sps} \isoto (S^{K_1} (\Ses))_{\Sps}
\]
between the formal completions of $S^K (\FS)$ and of $S^{K_1} (\Ses)$ 
along $\Sps$ compatible with the immersions $\isp$ and $i_1$. 
By analytical invariance \cite[Thm.~5.1]{W9},
$f$ induces an isomorphism
\[
\Mgm (\Sps, i_{\sigma , p}^! \, j_! \, \BZ) \isoto
\Mgm (\Sps, i_1^! j_{1 !} \,  \BZ) \; .
\]
Now apply Corollary~\ref{2C}.
\end{Proof}

Note \cite[Rem.~3.2]{W9} that the left hand side of the isomorphism from
Theorem~\ref{3main} lies in the full triangulated sub-category
$DM^{eff}_{gm} (E (P , \FX))$ 
of $DM^{eff}_-(E (P , \FX))$, and hence in the category $DM_{gm} (E (P , \FX))$ 
of geometrical motives over 
$E (P , \FX)$. This latter category is a rigid tensor category 
\cite[Thm.~4.3.7~1.\ and 2.]{VSF}. 
In particular, there exists an internal $Hom$
functor 
\[
\uHom: DM_{gm} (E (P , \FX)) \times DM_{gm} (E (P , \FX)) 
\longto DM_{gm} (E (P , \FX)) \; .
\]
Writing $M^* := \uHom (M, \BZ (0))$,
we thus have $M = (M^*)^*$ for all $M \in DM_{gm} (E (P , \FX))$. \\

If $K$ is neat, then $S^K$ is in $Sm /E (P , \FX)$.
Recall that in this case, 
the \emph{motive with compact support of $\Sps$ 
and with coefficients in $i_{\sigma , p}^* \, j_* \BZ$} is defined as
\[
\Mcgm (\Sps, i_{\sigma , p}^* \, j_* \BZ) :=
\Mgm (\Sps, i_{\sigma , p}^! \, j_! \, \BZ)^* (\dim S^K) [2 \dim S^K]
\]
\cite[Def.~7.1]{W9}. It lies in the full triangulated sub-category
$DM^{eff}_{gm} (E (P , \FX))$ of $DM_{gm} (E (P , \FX))$.

\begin{Cor} \label{3cor}
Assume that $(P, \FX)$ satisfies $(+)$, and that $K$ is neat. 
Then there is a canonical isomorphism 
\[
\Mcgm (\Sps, i_{\sigma , p}^* \, j_* \BZ) \isoto 
\Hom (\BZ (\sigma),\Mcgm (S^{K_1})) [- \dim \sigma]  
\]
in $DM^{eff}_{gm}(E (P , \FX))$.
\end{Cor}

\begin{Proof} 
Since the group $\BZ (\sigma)$ is canonically self-dual, Theorem~\ref{3main}
tells us that
$\Mcgm (\Sps, i_{\sigma , p}^* \, j_* \BZ)$ is canonically isomorphic to 
\[
\Hom (\BZ (\sigma),\Mgm (S^{K_1})^*(\dim S^K)) [2 \dim S^K - \dim \sigma] \; .  
\]
Note that $\dim S^K = \dim S^{K_1}$.
According to \cite[Thm.~4.3.7~3.]{VSF}, there is a canonical 
isomorphism
\[
\Mcgm (S^{K_1}) \isoto \Mgm (S^{K_1})^* (\dim S^{K_1})[2 \dim S^{K_1}] \; .
\]
This gives the desired isomorphism.
\end{Proof}

Theorem~\ref{3main} and Corollary~\ref{3cor} identify the motives
occurring as graded pieces in the co-localization \cite[Sect.~3]{W9}
resp.\ localization \cite[Sect.~7]{W9} filtration of the 
\emph{boundary motive} $\dMgm (S^K)$ introduced in \loccit ,
and associated to the stratification of a toroidal compactification.
When $S^K$ is pure,
the filtration associated to the Baily--Borel compactification
to be discussed in the next section is usually more efficient.


\bigskip

%
%

\section{The Baily--Borel compactification of pure Shimura varieties}
\label{4}



Assume that $(P, \FX) = (G, \FH)$ is pure. 
Fix an open compact subgroup $K \subset G (\BA_f)$.
Let $(S^K)^*$ denote the \emph{Baily--Borel compactification}
of $S^K$ \cite{BaBo,AMRT}. 
It comes equipped with a natural stratification into locally closed strata.
The unique open stratum is $S^K$. 
Any stratum diffe\-rent from the generic
one is obtained as follows: 
Fix a 
proper admissible parabolic subgroup $Q$ of $G$ with
associated normal subgroup $P_1$. Fix a
rational boundary
component $(P_1, \FX_1)$ of $(G, \FH)$,
and an element $g \in G(\BA_f)$.
Define
$K' := g \cdot K \cdot g^{-1}$, and
$K_1  :=  P_1 (\BA_f) \cap K'$.
Denote by $W_1$
the unipotent radical of $P_1$, and by 
\[
\pi: (P_1,\FX_1) \longonto (G_1,\FH_1) = (P_1,\FX_1)/W_1
\]
the quotient
of $(P_1 , \FX_1)$ by $W_1$ \cite[Prop.~2.9]{P1}.
From the proof
of \cite[Lemma~4.8]{P1}, it follows that $W_1$
equals the unipotent radical of $Q$. 
According to \cite[7.6]{P1}, there is
a canonical morphism
\[
\ig : \Spi := \Spi (G_1, \FH_1) \longto 
(S^K)^* 
\]
whose image is locally closed, and identical to the stratum in question. 
Define the following group (cmp.\ \cite[Sect.~1]{BW}):
\[
H_Q := \Stab_{Q (\BQ)} (\FH_1) \cap P_1 (\BA_f) \cdot K' \; .
\]
Here, the symbol $\Stab_{Q (\BQ)} (\FH_1)$ denotes
the subgroup of $Q (\BQ)$ stabilizing $\FH_1$.
$H_Q$ acts by analytic automorphisms on
$\FH_1 \times P_1 (\BA_f) / K_1$. Hence the group
$\Delta_1 := H_Q / P_1(\BQ)$ acts naturally on
\[
\Spi (\BC) =
P_1(\BQ) \backslash ( \FH_1 \times P_1 (\BA_f) / K_1 ) \; .
\]
This action is one by automorphisms of schemes over $E(G, \FH)$ 
\cite[Prop.~9.24]{P1}.
By \cite[6.3]{P1},
it factors through a finite quotient of $\Delta_1$.
The quotient by this action is precisely the image of $\ig$:
\[
\vcenter{\xymatrix@R-10pt{
\Spi \ar@{->>}[r] \ar@/_2pc/[rr]_{\ig}&
          S_1^K := \Delta_1 \backslash \Spi 
          \ar@{^{ (}->}[r]^-{\ig} & (S^K)^*
\\}}
\]
By abuse of notation, we denote by the same letter $\ig$ the inclusion of
the stratum $S_1^K$ into $(S^K)^*$.
Note that $\Delta_1$ is an
arithmetic subgroup of $Q/P_1 (\BQ)$, which is neat if $K$ is. By
\cite[11.1~(c)]{BS}, such a group is of
type $FL$, and hence satisfies condition (B) from Section~\ref{1b}.
Denote by 
\[
j: S^K \longinto (S^K)^* 
\]
the open immersion, by $S^{K_1}$ the Shimura variety
$S^{K_1} (P_1 , \FX_1)$, by $u_1$ the dimension of the Abelian unipotent
group $U_1$, and by $\BZ (U_1 (\BR))$ the group of orientations 
of the vector space of its $\BR$-valued points.
Recall the definition of the $\Delta_1$-equivariant motive
$M(S^{K_1}) \in \Delta_1 \text{-} DM^{eff}_-(E (G , \FH))$, and of the functor
\[
R \Gamma (\Delta_1, \argdot):  
\Delta_1 \text{-} DM^{eff}_-(E (G , \FH)) \longto DM^{eff}_-(E (G , \FH))
\] 
(Section~\ref{1b}). Here is 
our main result concerning strata in the Baily--Borel compactification:

\begin{Thm} \label{4main}
Assume that $(G, \FH)$ satisfies $(+)$, and that $K$ is neat. 
Then there is a canonical isomorphism 
\[
\Mgm (S_1^K, i_g^! \, j_! \, \BZ) \isoto
R \Gamma \left( \Delta_1 , \Hom (\BZ (U_1 (\BR)),\Mgm (S^{K_1})) \right) [u_1]  
\]
in $DM^{eff}_-(E (G , \FH))$.
\end{Thm}

\begin{Proof}
Choose a $K$-admissible cone decomposition $\FS$ satisfying the
conditions of \cite[(3.9)]{P2}. In particular, all cones occurring in
$\FS$ are smooth, and the decomposition is \emph{complete}.
Let us denote by $S^K (\FS) := S^K (G,\FH,\FS)$ the toroidal compactification
associated to $\FS$. It is a smooth projective scheme over $E (G , \FH)$.
The identity on
$S^K$ extends uniquely to a surjective proper morphism
\[
p = p_\FS: S^K (\FS) \longonto (S^K)^* \; .
\]
Denote by $Z'$ the inverse image under $p$ of the stratum $S_1^K$,
by $i_{\FS}$ the immersion of $Z'$ into $S^K (\FS)$, and by $j_{\FS}$ 
the open immersion of $S^K$. We apply \emph{invariance under abstract
blow-up} \cite[Thm.~4.1]{W9}, which tells us that 
\[
p: \Mgm (Z', i_{\FS}^! \, j_{\FS \, !} \, \BZ) \longto 
\Mgm (S_1^K, i_g^! \, j_! \, \BZ)
\]
is an isomorphism. In order to identify the left hand side,
we now set up data (1)--(6) as in Section~\ref{2}. 
The factorization of $\pi: (P_1,\FX_1) \to (G_1,\FH_1)$
corresponding to the weight filtration of the unipotent radical $W_1$
gives the following:
\[
\vcenter{\xymatrix@R-10pt{
(P_1,\FX_1) \ar@{->>}[r]^-{\pi_t} \ar@/_2pc/[rr]_\pi&
(P'_1,\FX'_1) := (P_1,\FX_1)/U_1 \ar@{->>}[r]^-{\pi_a} &
(G_1,\FH_1)
\\}}
\]
On the level of Shimura varieties, we get:
\[
\vcenter{\xymatrix@R-10pt{
X:= S^{K_1} (P_1,\FX_1) \ar@{->>}[r]^-{\pi_t} \ar@/_2pc/[rr]_\pi &
B \ar@{->>}[r]^-{\pi_a} &
\Spi = \Spi (G_1,\FH_1)
\\}}
\]
By \cite[3.12--3.22~(a)]{P1}, $\pi_a$ is in a natural way a torsor under an
Abelian scheme, while $\pi_t$ is a torsor under a split torus $T$
of dimension $u_1$.
Then \cite[(3.10)]{P2}:
$\FS$ and the choice of the stratum
$S^K_1$ induce a rational partial polyhedral decomposition of
$Y_*(T)_\BR$, equally denoted by $\FS$, and a non-empty proper subset
$\FT$ of $\FS$ satisfying conditions $(\alpha)$--$(\gamma)$
from Section~\ref{2}. In addition,
the subset $D$ of $Y_*(T)_\BR$ is contractible and open. As for the group,
we take $H := \Delta_1$. Its action on $B$, $X$ and $T$ is the one
induced from the natural action of $H_Q$ on the Shimura data involved
in the above factorization of $\pi$. It satisfies condition $(\delta)$.
Finally, we set $S := S^K$, $S_\FS := S^K (\FS)$, and $S_\FT := Z'$.
The isomorphisms
\[ 
H \backslash (X_\FS)_{X_\FT} \isoto (S_\FS)_{S_\FT}
\; , \; 
H \backslash (X_{\FS - \{ 0 \}})_{X_\FT} \isoto (S_{\FS - \{ 0 \}})_{S_\FT}
\; , \; 
H \backslash X_\FT \isoto S_\FT 
\]
are those of \cite[p.~224]{P2}.
Now apply Theorem~\ref{2main}~(c).

By passing to simultaneous refinements of two cone decompositions,
one sees that our isomorphism does not depend on the choice of $\FS$.
\end{Proof}

By \cite[Rem.~3.2]{W9}, the left hand side of the isomorphism from
Theorem~\ref{4main} lies in
$DM^{eff}_{gm} (E (G , \FH))$, 
and hence in the category $DM_{gm} (E (G , \FH))$ 
of geometrical motives over 
$E (G , \FH)$. Recall the functor 
\[
L \Lambda (\Delta_1,\argdot): \Delta_1 \text{-} DM^{eff}_-(E (G , \FH)) \longto 
DM^{eff}_-(E (G , \FH)) 
\]
from Section~\ref{1b}.

\begin{Cor} \label{4cor}
Assume that $(G, \FH)$ satisfies $(+)$, and that $K$ is neat. 
Then there is a canonical isomorphism 
\[
\Mcgm (S^K_1, i_g^* \, j_* \BZ) \isoto 
L \Lambda 
\left( \Delta_1 , \Hom (\BZ (U_1 (\BR)),\Mcgm (S^{K_1})) \right) [-u_1]  
\]
in $DM^{eff}_{gm}(E (G , \FH))$.
\end{Cor}

\begin{Proof}
By definition \cite[Def.~7.1]{W9}, the left hand side of the desired isomorphism
is dual to
\[
\Mgm (S^K_1, i_g^! \, j_! \, \BZ) (- \dim S^K) [- 2 \dim S^K] \; .
\]
Let us compute the dual of the right hand side, twisted by
$\dim S^K$ and shifted by $2 \dim S^K - u_1$.
In order to simplify notation, let us agree to drop the self-dual
group $\BZ (U_1 (\BR))$.
Recall that the object 
\[
L \Lambda ( \Delta_1 ,\Mcgm (S^{K_1}))
\]
of $DM^{eff}_-(E (G , \FH))$ is represented by the complex 
\[
\left( F_{\ast} \otimes_\BZ \uC_* (L^c(S^{K_1})) \right)_{\Delta_1} \; ,
\]
where $L^c(S^{K_1})$ denotes the Nisnevich sheaf with transfers
of quasi-finite correspondences, and where $F_\ast \to \BZ$ is
a bounded resolution of the trivial $\BZ \Delta_1$-module $\BZ$ by free
$\BZ \Delta_1$-modules. By adjunction,
\[
\uHom \left( L \Lambda ( \Delta_1 ,\Mcgm (S^{K_1})) ,
\BZ (\dim S^K) [2 \dim S^K] \right)
\]
is thus represented by 
\[
\Hom (F_{\ast}, z^*(\argdot,S^{K_1}))^{\Delta_1} \; ,
\]
for any complex $z^*(\argdot,S^{K_1})$  in
$C^-(\Delta_1 \text{-} Shv_{Nis}(SmCor(E (G , \FH))))$ represen\-ting
\[
\uHom \left( \Mcgm (S^{K_1}) , \BZ (\dim S^{K_1}) [2 \dim S^{K_1}] \right)
\]
(recall that $\dim S^K = \dim S^{K_1}$).
We claim that the class of $z^*(\argdot,S^{K_1})$
is canonically isomorphic to the $\Delta_1$-equivariant motive
$\Mgm (S^{K_1})$. Admitting this for the moment, we
apply Corollary~\ref{1bC} to get a canonical isomorphism
\[
\uHom \left( L \Lambda ( \Delta_1 ,\Mcgm (S^{K_1})) ,
\BZ (\dim S^K) [2 \dim S^K] \right) \isoto
R\Gamma( \Delta_1, \Mgm (S^{K_1})) \; .
\]
But by Theorem~\ref{4main}, the right hand side is isomorphic to
\[
\Mgm (S_1^K, i_g^! \, j_! \, \BZ) [ - u_1] \; .
\]
To conclude the proof, it remains to identify the object
\[
\uHom \left( \Mcgm (S^{K_1}) , \BZ (\dim S^{K_1}) [2 \dim S^{K_1}] \right)
\]
of $\Delta_1 \text{-} DM^{eff}_-(E (G , \FH))$.
One first constructs an isomorphism
\[
\uHom \left( \Mgm (S^{K_1}) , \BZ (\dim S^{K_1}) [2 \dim S^{K_1}] \right)
\isoto \Mcgm (S^{K_1})
\]
as in the proof of \cite[Thm.~4.3.7~3.]{VSF}. This requires
a $\Delta_1$-equivariant
version of duality for bivariant cycle cohomology, which is
provided by \cite[Thm.~3.21~(i) and its proof]{S}. Next, adjunction
gives a canonical morphism
\[
\Mgm (S^{K_1}) \longto
\uHom \left( \Mcgm (S^{K_1}) , \BZ (\dim S^{K_1}) [2 \dim S^{K_1}] \right)
\]
in $\Delta_1 \text{-} DM^{eff}_-(E (G , \FH))$. That it is an
isomorphism can be checked in
the category $DM^{eff}_-(E (G , \FH))$
\cite[Prop.~3.1]{S}. Now apply
\cite[Thm.~4.3.7~3.]{VSF}.
\end{Proof}

A more conceptual proof of Corollary~\ref{4cor} could be given
if the ca\-te\-gory
$\Delta_1 \text{-} DM_{gm} (E (G , \FH))$ were known to be a rigid tensor
category (see Remark~\ref{1bGa}),
in which the $\Delta_1$-equivariant motive is dual to the 
$\Delta_1$-equivariant motive with compact support (with the usual
twist and shift).\\

Theorem~\ref{4main} and Corollary~\ref{4cor} identify the motives
occurring as graded pieces in the co-localization \cite[Sect.~3]{W9}
resp.\ localization \cite[Sect.~7]{W9} filtration of the 
boundary motive $\dMgm (S^K)$. Here, the filtration is 
associated to the stratification of the Baily--Borel compactification. \\

The isomorphisms from Theorem~\ref{4main} and Corollary~\ref{4cor}
induce isomorphisms for \emph{motivic cohomology}
\[
H^i(X,\BZ(j)) := Hom_{DM_{gm}} (\Mgm (X), \BZ(j)[i])
\]
and for \emph{motivic cohomology with compact support}
\[
H_c^i(X,\BZ(j)) := Hom_{DM_{gm}} (\Mcgm (X), \BZ(j)[i]) \; ,
\]
for $i,j \in \BZ$, or more generally,
isomorphisms in any category equipped with a functor from
$DM^{eff}_{gm}(E (G , \FH))$. Let us spell out, for example,
what Corollary~\ref{4cor} implies for motivic cohomology:

\begin{Cor} \label{4C}
In the situation of Corollary~\ref{4cor}, there is a canonical isomorphism
of Abelian groups between
\[
Hom_{DM_{gm}(E (G , \FH))} (\Mcgm (S^K_1, i_g^* \, j_* \BZ), \BZ(j)[i])
\]
and
\[
H_c^{i+u_1} (S^{K_1},\Delta_1,\BZ(j)) \otimes_\BZ \BZ (U_1 (\BR)) \; ,
\]
for $i,j \in \BZ$.
Here, the $H_c^* (\argdot,\Delta_1,\BZ(j))$ denote
\emph{equivariant motivic cohomology with compact support},
i.e., the cohomology associated
to the composition
\[
R \Gamma (\Delta_1,\argdot)
\circ RHom_{DM^{eff}_-(E (G , \FH))}(\argdot, \BZ(j))
\circ \Mcgm
\]
mapping $\Delta_1 \text{-} Sch / E (G , \FH)$ to the derived category
of Abelian groups.
\end{Cor}

\begin{Rem} \label{4rema}
The analogous formula for \'etale cohomology
holds thanks to the existence of the \emph{\'etale realization functor}
\cite[3.3]{VSF}. Note that equi\-variant \'etale cohomology
\[ 
H_{c,et}^{i+u_1} (S^{K_1},\Delta_1,\BZ / n \BZ(j)) 
\otimes_\BZ \BZ (U_1 (\BR)) 
\]
equals
\[
H_{c,et}^{i+u_1} (\Spi,\Delta_1,R \pi_! \BZ / n \BZ(j)) 
\otimes_\BZ \BZ (U_1 (\BR)) \; .
\]
Using duality for the cohomology of
unipotent groups, and for $\pi_*$ and $\pi_!\, $, we get canonically
\[
H_{c,et}^{i+u_1} (S^{K_1},\Delta_1,\BZ / n \BZ(j)) 
\otimes_\BZ \BZ (U_1 (\BR)) \cong
H_{c,et}^{i} (\Spi,\Delta_1,R \pi_* \BZ / n \BZ(j)) \; .
\]
The resulting isomorphism
\[
H_{c,et}^{i} (S^K_1, i_g^* \, j_* \BZ /n \BZ(j)) \isoto
H_{c,et}^{i} (\Spi,\Delta_1,R \pi_* \BZ / n \BZ(j)) 
\]
should be identical to the one obtained by applying the functor $H_{c,et}^{i}$
to the comparison isomorphism from the
main result of \cite{P2} in the special case of Tate coefficients
$\BZ /n \BZ(j)$.
In the same way, but by applying the \'etale realization to 
the isomorphism in $DM^{eff}_{gm}(\bar{E})$ induced by the one
from Corollary~\ref{4cor} ($\bar{E}:=$ the algebraic closure of
$E (G , \FH)$), one gets an isomorphism
\[
H_{c,et}^{i} (S^K_{1,\bar{E}}, i_g^* \, j_* \BZ /n \BZ(j)) \isoto
H_{c,et}^{i} (S^{\pi (K_1)}_{\bar{E}},\Delta_1,R \pi_* \BZ / n \BZ(j)) 
\]
of modules under the absolute Galois group of $E (G , \FH)$.
Again, this isomorphism should be compared to the one deduced from
the main result of \cite{P2}. 
\end{Rem}

\begin{Rem} \label{4rem}
There is a motivic version of Theorem~\ref{4main} for more general
than constant coefficients $\BZ$. More precisely, one fixes mixed
Shimura data $(P,\FX)$ whose underlying pure data are $(G,\FH)$.
Choose an open compact subgroup $K^P \subset P(\BA_f)$ mapping onto $K$
under the projection $p$ from $(P,\FX)$ to $(G,\FH)$. One defines the motive 
\[
\Mgm (S_1^K, i_g^! \, j_! \, p_! \, \BZ)  
\]
by using a toroidal compactification of the mixed Shimura variety
$S^{K^P}(P,\FX)$ (one shows that the choice of compactification does not
matter). The general version of the isomorphism of Theorem~\ref{4main}
then takes the following form:
\[
\Mgm (S_1^K, i_g^! \, j_! \, p_! \, \BZ)  
\isoto
R \Gamma \left( \Delta_1^P , \Hom (\BZ (U_1^P (\BR)),\Mgm (S^{K_1^P})) \right) [u_1^P] 
\]
Here, $U_1^P$ is the ``weight $-2$'' part of the boundary component $(P_1^P,\FX_1^P)$ of 
$(P,\FX)$ lying over $(P_1,\FX_1)$, $u_1^P$ is its dimension, and 
$S^{K_1^P}$ is the Shimura variety associated to $(P_1^P,\FX_1^P)$. 
Finally, $\Delta_1^P$ is defined in exact analogy to $\Delta_1$, by
replacing $Q$, $P_1$, $K$, and $\FH_1$ by the respective groups and spaces
occurring in the description of the boundary component $(P_1^P,\FX_1^P)$.
Details will be published elsewhere.
\end{Rem}


\bigskip

%
%

\end{document}